\documentclass{amsart}

\usepackage{amsmath,amsthm,amssymb,amsfonts,psfig,graphicx,psfrag}

\numberwithin{equation}{section}
\theoremstyle{plain}
\newtheorem{theo+}           {Theorem}      [section]
\newtheorem{prop+}  [theo+]  {Proposition}
\newtheorem{lemm+}  [theo+]  {Lemma}
\newtheorem{cor+}  [theo+]  {Corollary}
\newenvironment{theorem}{\begin{theo+}}{\end{theo+}}

\newenvironment{corollary}{\begin{cor+}}{\end{cor+}}

\newcommand{\de}{\delta}
\newcommand{\la}{\lambda}

\newcommand{\De}{\Delta}
\begin{document}
\baselineskip 18pt
\larger[2]
\title[Elliptic hypergeometric series]
{Elliptic hypergeometric series\\ on root systems}
\author{Hjalmar Rosengren}
\address{Department of Mathematics\\ Chalmers University of Technology and 
G\"oteborg University\\SE-412~96 G\"oteborg, Sweden}
\email{hjalmar@math.chalmers.se}

\begin{abstract}
We derive a number of summation and transformation formulas
for elliptic hyper\-geo\-metric series on the root systems $A_n$, $C_n$
and $D_n$. In the special cases of classical and $q$-series,
our approach leads to new elementary proofs of the corresponding
identities. 
\end{abstract}
\maketitle   

\section{Introduction}
The study of hypergeometric series has, since the days of Gauss, been an 
active field of research with an increasing number of applications.
In a wide sense, this field includes a large number of different extensions of
Gauss's hypergeometric series,
such as the basic hypergeometric series or $q$-series first studied
by Cauchy and Heine \cite{gr}. During the last 20 years, $q$-series have, 
for good reasons, been enormously popular. Apart from the intrinsic
beauty of the subject, the main reasons seem to
be, on the one hand, applications of such series in number theory
and combinatorics, on the other hand, their relation to various
integrable models in mathematical physics and to the related algebraic
structures known as quantum groups. 

Recently, a new natural extension of hypergeometric series
was introduced, namely, the elliptic or modular hypergeometric
series of Frenkel and Turaev \cite{ft}. Subsequent papers on this topic
include  \cite{ds1,ds2,r,sp1,sp,sp3,sz,w}. 
 The motivation came
from  statistical mechanics; more precisely, they
may be used to express elliptic $6j$-symbols, which are elliptic solutions of 
the Yang--Baxter equation found
by Date et al.\ \cite{d}. Moreover, elliptic hypergeometric series are 
intimately connected to Ruij\-senaars's elliptic gamma function \cite{ru},
cf.~also \cite{fv},
which also arose in connection with integrable systems.
Integrable models involving elliptic
functions and the related elliptic quantum groups \cite{f} are of great
current interest. It is likely that their 
relations to elliptic hypergeometric series go far beyond
what is presently known.

From a mathematical point of view, elliptic hypergeometric series are
very natural objects. 
Recall that a series $\sum_k a_k$ is called hypergeometric if
$f(k)=a_{k+1}/a_k$ is a rational function of $k$, and a $q$-series if
$f$ is a rational function of $q^k$ for a fixed $q$. This may be compared
with Weierstrass's theorem, stating that a meromorphic function in $z$
satisfying an algebraic addition theorem is either rational, rational
in $q^z$ (we may then write $q=e^{ih}$ and call the function trigonometric)
or an elliptic function. Elliptic hypergeometric series form the top
level of the  hierarchy ``rational --- trigonometric
--- elliptic'' within the framework of hypergeometric series. Namely, for
elliptic hypergeometric series the term ratio $f$ is 
 an elliptic function. (It may seem natural to, more generally, allow $f$ to
be doubly quasi-periodic, but it appears that in interesting cases it is not
only elliptic but also has modular invariance; cf.~\cite{sp} for a detailed
discussion.)

Another natural extension of hypergeometric series is known as
hypergeometric series on root systems. 
Loosely speaking, these are multivariable series
such that the Weyl denominator of a classical root 
system is responsible for the coupling of the summation indices.  
During the last 25 years, this type of series
(classical and $q$-) has been studied  intensively by
Gustafson, Milne and many other researchers. The theory has 
found many  applications, such as to  plane partitions 
\cite{gek,kr}
and to representation of integers as sums of squares \cite{m5}.
Hypergeometric series on root systems  are also closely related to 
Macdonald polynomials and other multivariable 
 ortho\-gonal poly\-nomials \cite{vd,g1,kn,ra,st} and to Selberg-type
integrals and  Macdonald--Morris-type constant term 
identities for root systems \cite{vd,gu2,gb,gi,m1}.
(The cited references are a very small number of 
representative papers from a large and  growing literature.)

Hypergeometric series on root systems first appeared
(in the case of the root system $A_n$) 
in the work of Biedenharn, Holman and Louck \cite{hbl}, 
where they were used to express multiplicity-free $6j$-symbols of the
group $\mathrm{SU}(n)$. Thus, both elliptic hypergeometric series
and hypergeometric series on root systems were motivated by
mathematical physics, and even more precisely by different extensions
of  Racah's expression for the classical $6j$-symbol as
a ${}_4F_3$ hypergeometric sum.

The subject of the present work is the natural unification of these
two  fields: \emph{elliptic hypergeometric series on root
systems}. 

It should be pointed out that creating a theory of elliptic hypergeometric
 series on root systems is not  merely a matter of plugging  an extra parameter
into existing works.  The reason is that the theory
of (one- and multivariable) hypergeometric series is usually built up
as a ladder. That is, one starts with the simplest results such as
the binomial theorem and Gauss's ${}_2F_1$ sum, and use these to successively
build more complicated identities. A natural culmination point is given 
by Jackson's
${}_8W_7$ summation and Bailey's ${}_{10}W_9$ transformation (or their
multivariable analogues). This is a natural approach which follows the 
historical development, but it is doomed to fail for elliptic hypergeometric
series (which is surely one reason why it took so long before 
these were discovered). 
The problem is that, for all known elliptic hypergeometric series identities,
the simplest non-trivial case is  the addition formula
\begin{equation}\label{aa}\begin{split}&
\theta(x+z)\,\theta(x-z)\,\theta(y+w)\,\theta(y-w)\\
&\ =\theta(x+y)\,\theta(x-y)\,\theta(z+w)\,\theta(z-w)+\theta(x+w)\,\theta(x-w)
\,\theta(y+z)\,\theta(y-z)
\end{split}\end{equation}
for a suitably defined theta function $\theta$ (this is not an addition
theorem in the sense of Weierstrass).
In our figurative ladder, the lower steps correspond to lower level
addition formulas such as
$$\sin(x+z)\sin(x-z)=\sin(x+y)\sin(x-y)+\sin(y+z)\sin(y-z), $$
which do not admit elliptic analogues. This means that the
theory of elliptic hypergeometric series looks like a ladder with all
but the top rungs missing. Thus, instead of climbing the ladder one must
find a way to start at the top. This may be viewed as 
additional motivation for studying elliptic hypergeometric series, since
such an alternative approach should add something to our understanding
also of the rational and trigonometric levels. 

The present paper is not the first
 to deal with elliptic hypergeometric series on root systems.
In \cite{w}, Warnaar proved one multivariable elliptic analogue of Jackson's
${}_8W_7$ summation formula (cf.~\eqref{wj} below) 
and conjectured another one. We will refer to these as Warnaar's first
and second summation.
Van Diejen and Spiridonov \cite{ds1} found that the second
summation may be derived from a certain conjectured elliptic Selberg 
integral. In \cite{r}, we gave an inductive proof of Warnaar's
second summation, using a very special case of his first  summation in 
the inductive step. It should be said that although both
Warnaar sums live on root systems (to be precise, on $C_n$), 
they are somewhat different in nature
from those going back to the seminal paper \cite{hbl}. 
In the terminology of \cite{ds2}, the first sum is of
 ``Schlosser-type'' and the second one of ``Aomoto--It\^o--Macdonald-type'',
as opposed to the Gustafson--Milne-type series which are our main concern
here. In the latter direction, 
van Diejen and Spiridonov \cite{ds2} stated a third
elliptic $C_n$ analogue of Jackson's 
summation  formula (cf.~Theorem \ref{cjt} below), 
which in the trigonometric case is due independently to Denis and
Gustafson \cite{dg} and  to Milne and Lilly \cite{ml}. However,
the identity was only proved under the condition that a certain
Selberg-type integral evaluation conjectured in \cite{ds} holds.

The purpose of the present paper is to give elliptic analogues of
a a number of known identities for Gustafson--Milne-type series
on the root systems $A_n$, $C_n$ and $D_n$. 
We have focused on deriving some of the most fundamental 
and important identities --- analogues of Jackson's ${}_8W_7$
summation and Bailey's ${}_{10}W_9$ transformation, although it seems
likely that more identities (for instance, 
  quadratic and bibasic summation formulas) 
may be derived by our methods. In the one-variable case, our identities
reduce to either the elliptic Jackson summation
 or the elliptic Bailey transformation of Frenkel and Turaev \cite{ft}.

To obtain elliptic Jackson summations on
 the root systems $A_n$ and $D_n$ we use an inductive method similar
to our proof of Warnaar's second summation in \cite{r}. Instead of applying
 Warnaar's first summation in the induction step, we
use elliptic analogues of elementary partial fraction expansions. 
In our opinion, the resulting proofs are simpler than those previously
known for classical and $q$-series,  even for very degenerate cases of the
 identities.
For the root system $C_n$, we cannot use this approach. Instead, we
observe that the elliptic $C_n$ Jackson summation we want to prove
(the same one that was conjectured in \cite{ds2}) may be obtained, in a rather
non-obvious way, as a special case of Warnaar's first summation. 
This gives an unexpected link with hypergeometric series connected with
determinant evaluations (the Schlosser-type series alluded to above).

The plan of the paper is as follows. In Section 2, we give a first illustration
of our new approach to $A_n$  and $D_n$ series, by showing how
it works in a simple  case, Milne's $A_n$ analogue of the
terminating ${}_6W_5$ summation. 
Sections 3 and  4 contain notation and preliminaries
for theta functions and elliptic hypergeometric series. Our main results
are contained in Sections 5,  6 and  7, where we derive elliptic Jackson 
summations on the root systems $A_n$, $D_n$ and $C_n$, respectively. In 
Section 8 we apply
these summations to derive a number of elliptic Bailey transformations on
these root systems. Since the step from
${}_8W_7$ sums to ${}_{10}W_9$ transformations works exactly as in the
rational and trigonometric case,  we do not provide much details in this
final section.

{\bf Remark and acknowledgement:} While I was finishing this paper, 
Vyacheslav Spiridonov informed me that he has independently conjectured
the summation formulas given here as  Corollary \ref{cr} and Corollary 
\ref{drs}
and, more importantly, verified that both sides are Jacobi modular 
forms in the sense of \cite{ez}. This is another strong indication 
that these identities are ``natural''. 
I would like to thank Prof.\ Spiridonov for this information, as 
well as for providing me with copies of the papers \cite{sp,sp3}.

\section{Milne's fundamental theorem}

To illustrate our method, we start with  one of the simplest 
identities for hyper\-geometric series on root systems, 
Milne's fundamental theorem of $\mathrm{U}(n)$ series.
\emph{In this section only} we use the notation
\begin{equation}\label{tp}(a)_k=(1-a)(1-aq)\dotsm(1-aq^{k-1}), \end{equation}
where $q$ is a fixed parameter suppressed from the notation. We also write,
again in this section only,
$$\De(z)=\prod_{1\leq j<k\leq n}(z_j-z_k)$$
for the Vandermonde determinant, which is the Weyl denominator for the
root system $A_n$. Basic hypergeometric series on  $A_n$ 
contain the factor
$$\frac{\De(zq^y)}{\De(z)}=\prod_{1\leq j<k\leq n}\frac{z_jq^{y_j}-z_kq^{y_k}}
{z_j-z_k},$$
where the $z_j$ are fixed parameters and the $y_j$ summation indices. 

We can now state Milne's fundamental theorem:
\begin{equation}\label{ft}
\sum_{\substack{y_1,\dots,y_n\geq 0\\y_1+\dotsm+y_n=N}}
\frac{\De(zq^y)}{\De(z)}\prod_{j,k=1}^n\frac{(a_jz_k/z_j)_{y_k}}
{(qz_k/z_j)_{y_k}}= \frac{(a_1\dotsm a_n)_N}{(q)_N}.\end{equation}
When $n=2$ this is Rogers's terminating ${}_6W_5$ summation formula
\cite[Equation (II.20)]{gr}. On the other hand,  multiplying \eqref{ft}
with $t^N$ and then summing over non-negative integers $N$ gives an
extension of Cauchy's $q$-binomial theorem \cite[Equation (II.3)]{gr}.
So one may also view \eqref{ft} as a refined multivariable $q$-binomial
theorem. 

 The identity \eqref{ft} was 
first obtained in \cite{m1}, where it was 
used to derive the Macdonald identities for the affine Lie algebra 
$A^{(1)}_n$. In Milne's approach to hypergeometric series on $A_n$ or 
$\mathrm{U}(n)$
(cf.~\cite{ms} for a recent survey), \eqref{ft} is the starting point
from which all other results are built; this motivates the name
``fundamental theorem of $\mathrm{U}(n)$ series'' \cite{b}. 
By contrast, in the present approach
 \eqref{ft} appears as a degenerate case of the elliptic Jackson summation
given in Theorem \ref{taj} below, which we will prove by precisely 
the same method. We only give \eqref{ft} a separate treatment in order to
illustrate our method in the simplest possible situation.

Milne's proof of \eqref{ft} is based on showing that both sides 
satisfy the same difference equation in the variables $a_k$ and        
$z_k$. An important tool is an identity which we write here as
\begin{equation}\label{pf}\sum_{k=1}^n
\frac{\prod_{j=1}^n(b_j-a_k)}{a_k\prod_{j\neq k}(a_j-a_k)}=
\frac{b_1\dotsm b_n}{a_1\dotsm a_n}-1\end{equation}
(throughout the paper we tacitly assume that all expressions are well-defined,
in this case that the $a_j$ are different and non-zero).
Our proof is also based on \eqref{pf}, but instead of difference
operators we use a simple inductive argument. 

The easiest way to prove
\eqref{pf} is to rewrite it as a  partial fraction expansion,
replacing $a_j$, $b_j$, by $a_j-t$, $b_j-t$: 
$$\prod_{j=1}^n\frac{b_j-t}{a_j-t}=1+\sum_{k=1}^n
\frac{\prod_{j=1}^n(b_j-a_k)}{(a_k-t)\prod_{j\neq k}(a_j-a_k)}.$$
The proof is then  easy: the 
 existence of the expansion is immediate by induction on
$n$, and to compute the coefficients one may multiply with $a_k-t$
and plug in $t=a_k$. This is taught at Swedish technical universities
under the name ``handp\aa l\"aggning'' (laying on hands) \cite{ap}.

We prove \eqref{ft} by induction on $N$, starting with  $N=1$.
In that case we have $y_i=\de_{ik}$ for some $k$. Using $k$ as
summation index, the left-hand side of \eqref{ft} may be 
rewritten and summed using \eqref{pf}:
$$\sum_{k=1}^n\prod_{j\neq k}\frac{z_j-qz_k}{z_j-z_k}\prod_{j=1}^n
\frac{z_j-a_jz_k}{z_j-qz_k}=\frac{a_1\dotsm a_n}{1-q}\sum_{k=1}^n
\frac{\prod_{j=1}^n(z_j/a_j-z_k)}{z_k\prod_{j\neq k}(z_j-z_k)}
=\frac{1-a_1\dotsm a_n}{1-q}.$$
 In fact, we see that the case $N=1$ of \eqref{ft}
is  equivalent to \eqref{pf}.  

Next we assume that \eqref{ft} holds for a fixed value of $N$. We may 
then write the right-hand side, with $N$ replaced by $N+1$, as
\begin{equation*}\begin{split}R&= \frac{(a_1\dotsm a_n)_{N+1}}{(q)_{N+1}}=
\frac{1-q^Na_1\dotsm a_n}{1-q^{N+1}}\frac{(a_1\dotsm a_n)_{N}}{(q)_{N}}\\
&=\frac{1-q^Na_1\dotsm a_n}{1-q^{N+1}}
\sum_{\substack{y_1,\dots,y_n\geq 0\\y_1+\dotsm+y_n=N}}
\frac{\De(zq^y)}{\De(z)}\prod_{j,k=1}^n\frac{(a_jz_k/z_j)_{y_k}}
{(qz_k/z_j)_{y_k}}.\end{split}\end{equation*}
We now  split the terms 
 using the case $N=1$ of \eqref{ft}, with  
$a_k$ replaced by $a_kq^{y_k}$ and $z_k$ by $z_kq^{y_k}$. Explicitly, we have
$$\frac{1-q^Na_1\dotsm a_n}{1-q}=
\sum_{\substack{x_1,\dots,x_n\geq 0\\x_1+\dotsm+x_n=1}}
\frac{\De(zq^{y+x})}{\De(zq^y)}\prod_{j,k=1}^n\frac{(q^{y_k}a_jz_k/z_j)_{x_k}}
{(q^{1+y_k-y_j}z_k/z_j)_{x_k}}, $$
which gives
\begin{equation}\label{r}\begin{split}R&=\frac{1-q}{1-q^{N+1}}
\sum_{\substack{y_1,\dots,y_n\geq 0\\y_1+\dotsm+y_n=N}}
\sum_{\substack{x_1,\dots,x_n\geq 0\\x_1+\dotsm+x_n=1}}
\frac{\De(zq^{y+x})}{\De(z)}\prod_{j,k=1}^n\frac{(a_jz_k/z_j)_{y_k+x_k}}
{(qz_k/z_j)_{y_k}(q^{1+y_k-y_j}z_k/z_j)_{x_k}}\\
&=\frac{1-q}{1-q^{N+1}}\sum_{\substack{y_1,\dots,y_n\geq 0\\
y_1+\dotsm+y_n=N+1}}\frac{\De(zq^{y})}{\De(z)}
\prod_{j,k=1}^n(a_jz_k/z_j)_{y_k}\\
&\quad\times\sum_{\substack{0\leq x_k\leq y_k\\x_1+\dotsm+x_n=1}}
\prod_{j,k=1}^n\frac{1}
{(qz_k/z_j)_{y_k-x_k}(q^{1+y_k-y_j-x_k+x_j}z_k/z_j)_{x_k}},
\end{split}\end{equation}
where we replaced $y$ by $y-x$. 

We will identify the sum in $x$ as a special case of \eqref{pf}.
Since $x_k\in\{0,1\}$, we have
\begin{equation}\label{e1}\frac{1}{(qz_k/z_j)_{y_k-x_k}}=
\frac{(q^{1+y_k-x_k}z_k/z_j)_{x_k}}{(qz_k/z_j)_{y_k}}
=\frac{(q^{y_k}z_k/z_j)_{x_k}}{(qz_k/z_j)_{y_k}}.\end{equation}
Writing the factor in that form makes the conditions $x_k\leq y_k$
superfluous, since $(q^{y_k})_{x_k}$
vanishes if $y_k=0$ and $x_k=1$.
Moreover, if $j\neq k$ the factor $(q^{1+y_k-y_j-x_k+x_j}z_k/z_j)_{x_k}$
equals $1$ unless $x_k=1$ and $x_j=0$, which gives 
\begin{equation}\label{e2}\prod_{j,k=1}^n(q^{1+y_k-y_j-x_k+x_j}z_k/z_j)_{x_k}
=(1-q)\prod_{j\neq k}
(q^{y_k-y_j}z_k/z_j)_{x_k}.\end{equation}
We may now compute the sum in $x$ using \eqref{pf}:
\begin{equation*}\begin{split}&
\frac{1}{1-q}\prod_{j,k=1}^n\frac{1}{(qz_k/z_j)_{y_k}}
\sum_{\substack{x_1,\dots,x_n\geq 0\\x_1+\dotsm+x_n=1}}
\frac{\prod_{j,k=1}^n(q^{y_k}z_k/z_j)_{x_k}}
{\prod_{j\neq k}(q^{y_k-y_j}z_k/z_j)_{x_k}}\\ 
&\quad=\frac{q^{N+1}}{1-q}\prod_{j,k=1}^n\frac{1}{(qz_k/z_j)_{y_k}}
\sum_{k=1}^n\frac{\prod_{j=1}^n(z_j-q^{y_k}z_k)}
{q^{y_k}z_k\prod_{j\neq k}(q^{y_j}z_j-q^{y_k}z_k)}\\
&\quad=\frac{1-q^{N+1}}{1-q}\prod_{j,k=1}^n\frac{1}{(qz_k/z_j)_{y_k}}.
\end{split}\end{equation*}
Plugging this into \eqref{r} gives
$$R=\sum_{\substack{y_1,\dots,y_n\geq 0\\y_1+\dotsm+y_n=N+1}}
\frac{\De(zq^y)}{\De(z)}\prod_{j,k=1}^n\frac{(a_jz_k/z_j)_{y_k}}
{(qz_k/z_j)_{y_k}},$$
which is the left-hand side of \eqref{ft}  with $N$ replaced by $N+1$.
This completes our  proof of Milne's fundamental theorem.

\section{Notation and preliminaries}
Elliptic hypergeometric series may be built from the theta function
$$\theta(x)=\prod_{j=0}^\infty(1-p^jx)(1-p^{j+1}/x),\qquad |p|<1,$$
where $p$ is a constant which is fixed throughout the paper.
It satisfies the inversion formula 
\begin{equation}\label{ti}\theta(x)=-x\,\theta(1/x),\end{equation}
the quasi-periodicity
\begin{equation}\label{tqp}\theta(px)=-\frac 1 x\,\theta(x)\end{equation}
and the ``addition formula''
\begin{equation}\label{ra}\begin{split}&
\theta(xz)\,\theta(x/z)\,\theta(yw)\,\theta(y/w)\\
&\quad=\frac yz\, \theta(xy)\,\theta(x/y)\,\theta(zw)\,\theta(z/w) 
+\theta(xw)\,\theta(x/w)\,\theta(yz)\,\theta(y/z)\end{split}\end{equation}
(cf.~\cite[Exercise 2.16 and Exercise 5.21]{gr} for proofs within the 
framework of basic hypergeometric series).
Equivalently, $x\mapsto q^{-x/2}\theta(q^x)$, with $q$ fixed, satisfies
 \eqref{aa}.  All identities for elliptic
hypergeometric series obtained in this paper may be viewed as  generalizations
of \eqref{aa}. 

We define elliptic Pochhammer symbols by
$$(a)_k=\theta(a)\,\theta(aq)\dotsm\theta(aq^{k-1}).$$
Note that they depend on two parameters $p$, $q$, which will be suppressed
from the notation. When
 $p=0$, $\theta(x)=1-x$ and we recover the symbols \eqref{tp}
used in the previous section. The elliptic symbols satisfy similar 
elementary identities as in the case $p=0$, such as
\begin{equation}\label{ap}(a)_{n+k}=(a)_n(aq^n)_k,\end{equation}
\begin{equation}\label{ep}(a)_{n-k}=(-1)^kq^{\binom k2}(q^{1-n}/a)^k
\frac{(a)_n}{(q^{1-n}/a)_k},\end{equation}
\begin{equation}\label{ip}(a)_n=(-1)^n q^{\binom n2}a^n(q^{1-n}/a)_n;
\end{equation}
these will be used repeatedly and usually without comment.
We  also need the quasi-periodicity
\begin{equation}\label{qp}(pa)_k=(-1)^kq^{-\binom k2}a^{-k}(a)_k,
\end{equation}
which follows from \eqref{tqp}.

We introduce the elliptic Vandermonde determinant
$$\De(z)=\prod_{1\leq j<k\leq n} z_j\,\theta(z_k/z_j),$$
which appears in elliptic hypergeometric series on $A_n$,
and, together with additional double products, 
in the series on $C_n$ and $D_n$.

The following identity will be useful
(cf.~\cite[Lemma 4.3]{m} for a  more general identity
in the  case $p=0$):
\begin{equation}\label{di}(-1)^{|y|}q^{|y|+\binom{|y|}{2}}\frac{\De(zq^y)}
{\De(z)}\prod_{j,k=1}^n\frac{(q^{-y_j}z_k/z_j)_{y_k}}{(qz_k/z_j)_{y_k}}=1.
\end{equation}
To prove this, one may rewrite the double product as
$$\prod_{j,k=1}^n\frac{(q^{-y_j}z_k/z_j)_{y_k}}{(qz_k/z_j)_{y_k}}
=\prod_{k=1}^n\frac{(q^{-y_k})_{y_k}}{(q)_{y_k}}\prod_{1\leq j<k\leq n} 
\frac{(q^{-y_j}z_k/z_j)_{y_k}(q^{-y_k}z_j/z_k)_{y_j}}{(qz_k/z_j)_{y_k}
(qz_j/z_k)_{y_j}},$$
which, by  \eqref{ip}, equals
\begin{equation*}\begin{split}& \prod_{k=1}^n(-1)^{y_k}q^{\binom {y_k}2-y_k^2}
\prod_{1\leq j<k\leq n}q^{-(y_k+1)y_j} 
\frac{(q^{-y_j}z_k/z_j)_{y_k}(q^{1-y_j+y_k}z_k/z_j)_{y_k}}{(qz_k/z_j)_{y_k}
(q^{-y_j}z_k/z_j)_{y_k}}\\
&\quad= (-1)^{|y|}q^{\sum_k\binom {y_k}2-y_k^2}\,q^{-\sum_{j<k}y_ky_j}
\prod_{1\leq j<k\leq n}q^{-y_j}\frac{\theta(z_k/z_j)}
{\theta(q^{y_k-y_j}z_k/z_j)}\\
&\quad=(-1)^{|y|}
q^{-\binom{|y|}{2}-|y|}\frac{\De(z)}{\De(zq^y)}.
\end{split}\end{equation*}

Finally, we will often use the shorthand notation
$$\theta(a_1,\dots,a_n)=\theta(a_1)\dotsm\theta(a_n),$$
$$(a_1,\dots,a_n)_k=(a_1)_k\dotsm(a_n)_k.$$

\section{Elliptic partial fraction expansions}

The starting point for our investigation of elliptic hypergeometric series 
on $A_n$ will be
the following elliptic generalization of \eqref{pf}: 
\begin{equation}\label{ww}\sum_{k=1}^n\frac{\prod_{j=1}^n\theta(a_k/b_j)}
{\prod_{j\neq k}\theta(a_k/a_j)}=0,\qquad a_1\dotsm a_n=b_1\dotsm b_n,
\end{equation}
which is valid as long as both sides are well-defined.
We do not know who first wrote down this identity.
It is given as an exercise in the classic textbook of
Whittaker and Watson \cite[p.~451]{ww}. 

To see that \eqref{ww} generalizes \eqref{pf}, replace $n$ by 
$n+1$ and write $a_{n+1}=t$.
Isolating the term with $k=n+1$, \eqref{ww} may be written as
 \begin{equation}\label{epf}
\sum_{k=1}^n\frac{\prod_{j=1}^{n+1}\theta(a_k/b_j)}
{\theta (a_k/t)\prod_{j\neq k}\theta(a_k/a_j)}=
-\frac{\prod_{j=1}^{n+1}\theta(t/b_j)}
{\prod_{j=1}^n\theta(t/a_j)},\qquad b_1\dotsm b_{n+1}=a_1\dotsm a_nt
\end{equation}
(note that one may recover \eqref{ww} also by letting $t=b_{n+1}$ in
\eqref{epf}).
If we put $p=0$ and then let $t,\,b_{n+1}\rightarrow 0$ 
with $t/b_{n+1}=b_1\dotsm b_n/a_1\dotsm a_n$ fixed
we recover \eqref{pf}. We also observe that, by \eqref{ti}, 
\eqref{epf} may alternatively be written as
 \begin{equation}\label{epa}
\sum_{k=1}^n\frac{\prod_{j=1}^{n+1}\theta(a_k/b_j)}
{\theta (a_k/t)\prod_{j\neq k}\theta(a_k/a_j)}=
\frac{\prod_{j=1}^{n+1}\theta(b_j/t)}
{\prod_{j=1}^n\theta(a_j/t)},\qquad b_1\dotsm b_{n+1}=a_1\dotsm a_nt.
\end{equation}

 A simple  proof of \eqref{epf} arises from the interpretation
 as a generalized partial fraction expansion. 
Namely, for $n=2$, \eqref{epf} is equivalent to the addition
formula \eqref{ra}. It then follows immediately by induction on $n$ that there
exists an expansion
\begin{equation}\label{pfa}\frac{\prod_{j=1}^n\theta(t/b_j)}
{\prod_{j=1}^n\theta(t/a_j)}=
\sum_{k=1}^n C_k\,\frac{\theta(a_kb_1\dotsm b_n/a_1\dotsm a_nt)}
{\theta (a_k/t)},\end{equation}
where the coefficients  are easily computed using ``handp\aa l\"aggning''. 
That is, if one 
multiplies \eqref{pfa} with  $\theta(t/a_k)=-\theta(a_k/t)t/a_k$
and let $t=a_k$ one obtains
$$C_k=-\frac{\prod_{j=1}^n\theta(a_k/b_j)}
{\theta(b_1\dotsm b_n/a_1\dotsm a_n)\prod_{j\neq k}\theta(a_k/a_j)}.$$
Plugging this into \eqref{pfa} and writing
 $b_{n+1}=a_1\dotsm a_nt/b_1\dotsm b_{n}$ yields \eqref{epf}.

For the root system $D_n$ we need in addition to \eqref{ww} the  identity
\begin{equation}\label{pfd}
\sum_{k=1}^n\frac{\prod_{j=1}^{n-1}\theta(a_kb_j,a_k/b_j)}
{\theta(ta_k,t/a_k)\prod_{j\neq k}\theta(a_ka_j,a_k/a_j)}
=\frac{\prod_{j=1}^{n-1}\theta(tb_j,t/b_j)}{\prod_{j=1}^n\theta(ta_j,t/a_j)},
\end{equation}
or, equivalently, 
\begin{equation}\label{dea}
\sum_{k=1}^n\frac{a_k\prod_{j=1}^{n-2}\theta(a_kb_j,a_k/b_j)}
{\prod_{j\neq k}\theta(a_ka_j,a_k/a_j)}=0.\end{equation}
The earliest occurrence of  \eqref{pfd} that we have found
is as \cite[Lemma 4.14]{gu}, where an analytic proof is given. 
Again, the interpretation  as 
 a generalized partial fraction expansion  makes an algebraic proof easy: 
for $n=2$ \eqref{pfd} is 
equivalent to \eqref{ra}, the existence of such an expansion 
follows by induction on $n$ and the coefficients may be obtained by
multiplying with $\theta(ta_k,t/a_k)$ and then letting $t=a_k$.

Finally, we point out some occurrences of 
the identities \eqref{ww} and \eqref{dea} in the related literature.
In \cite{ds}, they are given as Lemma A.2 and Lemma A.1, and, written
in the alternative 
form \eqref{epf} and \eqref{pfd}, as Proposition A.4 and Proposition
A.3. Implicitly, \eqref{dea} also appears in 
 \cite{w}. Namely, in \cite[Lemma 3.2]{w}  a certain matrix inversion, 
say $AB=I$, is obtained as a special case of a more general result. 
The inverse relation $BA=I$ is 
easily seen to be equivalent to \eqref{dea}. 

\section{An elliptic $A_n$ Jackson summation}

The first main result of the paper is the following Jackson-type summation
formula for an elliptic hypergeometric series on the root system $A_n$.
When $p=0$ we recover (with a new, elementary proof) Milne's $A_n$ 
Jackson summation \cite[Theorem 6.17]{mr}. To see this, one should
rewrite our identity as in Corollary \ref{aaj} below and Milne's
identity as in \cite[Theorem A.5]{mn}. 
Note also that \eqref{ft}
may be obtained from the case $p=0$ of Theorem \ref{taj} by
letting $b,\,a_{n+1}\rightarrow 0$ with $b/a_{n+1}$ fixed.

\begin{theorem}\label{taj}
If $b=a_1\dotsm a_{n+1}z_1\dotsm z_n$, the following identity holds:
$$\sum_{\substack{y_1+\dotsm+y_n=N\\y_1,\dots,y_n\geq 0}}
\frac{\De(zq^y)}{\De(z)}\prod_{k=1}^n\frac{\prod_{j=1}^{n+1}(a_jz_k)_{y_k}}
{(bz_k)_{y_k}\prod_{j=1}^n(qz_k/z_j)_{y_k}}=\frac{(b/a_1,\dots,b/a_{n+1})_N}
{(q,bz_1,\dots,bz_n)_N}.$$
\end{theorem}

\begin{proof}
We will use induction on $N$. If $N=1$, we have $y_i=\de_{ik}$ for
some $k$; using
$k$ as summation variable yields the identity
$$\sum_{k=1}^n\frac{\prod_{j=1}^{n+1}\theta(a_jz_k)}{\theta(bz_k)
\prod_{j\neq k}\theta(z_k/z_j)}=\frac{\prod_{j=1}^{n+1}\theta(b/a_j)}
{\prod_{j=1}^n\theta(bz_j)},\qquad b=a_1\dotsm a_{n+1}z_1\dotsm z_n, $$
which is equivalent to \eqref{epa}.

Assume that the theorem holds for a fixed $N$, and let 
$$R=\frac{(b/a_1,\dots,b/a_{n+1})_{N+1}}
{(q,bz_1,\dots,bz_n)_{N+1}}.$$
Imitating the proof of \eqref{ft} given above, we rewrite $R$
using the induction hypothesis, and then expand each term using
the case $N=1$, $z_k\mapsto q^{y_k}z_k$, $b\mapsto q^Nb$
of the theorem.  Explicitly, this gives
\begin{equation*}\begin{split}
R&=\frac{\prod_{k=1}^{n+1}\theta(q^Nb/a_k)}
{\theta(q^{N+1})\prod_{k=1}^n\theta(q^{N}bz_k)} 
\frac{(b/a_1,\dots,b/a_{n+1})_N}
{(q,bz_1,\dots,bz_n)_N}\\
&=\frac{\prod_{k=1}^{n+1}\theta(q^Nb/a_k)}
{\theta(q^{N+1})\prod_{k=1}^n\theta(q^{N}bz_k)}
\sum_{\substack{y_1+\dotsm+y_n=N\\y_1,\dots,y_n\geq 0}}
\frac{\De(zq^y)}{\De(z)}\prod_{k=1}^n\frac{\prod_{j=1}^{n+1}(a_jz_k)_{y_k}}
{(bz_k)_{y_k}\prod_{j=1}^n(qz_k/z_j)_{y_k}}\\
&=\frac{\theta(q)}{\theta(q^{N+1})}
\sum_{\substack{y_1+\dotsm+y_n=N\\y_1,\dots,y_n\geq 0}}
\Bigg(\frac{\De(zq^y)}{\De(z)}\prod_{k=1}^n\frac{\theta(q^{N+y_k}bz_k)
\prod_{j=1}^{n+1}(a_jz_k)_{y_k}}
{\theta(q^{N}bz_k)\,(bz_k)_{y_k}\prod_{j=1}^n(qz_k/z_j)_{y_k}}\\
&\quad\times
\sum_{\substack{x_1+\dotsm+x_n=1\\x_1,\dots,x_n\geq 0}}
\frac{\De(zq^{y+x})}{\De(zq^y)}
\prod_{k=1}^n\frac{\prod_{j=1}^{n+1}(q^{y_k}a_jz_k)_{x_k}}
{(q^{N+y_k}bz_k)_{x_k}\prod_{j=1}^n(q^{1+y_k-y_j}z_k/z_j)_{x_k}}\Bigg).
\end{split}\end{equation*}

Next we replace $y$ by $y-x$ in the summation. 
Since $x_k\in\{0,1\}$, we may write
$$\frac{\theta(q^{N+y_k-x_k}bz_k)}{(q^{N+y_k-x_k}bz_k)_{x_k}}=
\frac{\theta(q^{N+y_k}bz_k)}{(q^{N+y_k}bz_k)_{x_k}},$$
$$\frac{1}{(bz_k)_{y_k-x_k}}=\frac{(q^{y_k-1}bz_k)_{x_k}}{(bz_k)_{y_k}}.$$
It is clear that \eqref{e1} holds also in the elliptic case, and
that \eqref{e2} has the elliptic analogue
$$\prod_{j,k=1}^n(q^{1+y_k-y_j-x_k+x_j}z_k/z_j)_{x_k}
=\theta(q)\prod_{j\neq k}(q^{y_k-y_j}z_k/z_j)_{x_k}.$$
This gives the expression
\begin{equation}\label{re}\begin{split}
R&=\frac{1}
{\theta(q^{N+1})}
\sum_{\substack{y_1+\dotsm+y_n=N+1\\y_1,\dots,y_n\geq 0}}
\Bigg(\frac{\De(zq^y)}{\De(z)}\prod_{k=1}^n
\frac{\theta(q^{N+y_k}bz_k)\prod_{j=1}^{n+1}(a_jz_k)_{y_k}}
{\theta(q^Nbz_k)\,(bz_k)_{y_k}\prod_{j=1}^n(qz_k/z_j)_{y_k}}\\
&\quad\times\sum_{\substack{x_1+\dotsm+x_n=1\\x_1,\dots,x_n\geq 0}}
\prod_{k=1}^n\frac{(q^{y_k-1}bz_k)_{x_k}\prod_{j=1}^n(q^{y_k}z_k/z_j)_{x_k}}
{(q^{N+y_k}bz_k)_{x_k}\prod_{j\neq k}
(q^{y_k-y_j}z_k/z_j)_{x_k}}\Bigg).
\end{split}\end{equation}
  (Formally, we have introduced extra  terms corresponding to $x_k=1$ and 
$y_k=0$, however, these  all
vanish in view of the factor $(q^{y_k}z_k/z_j)_{x_k}$.)
The sum in $x$ may be rewritten and summed using \eqref{epa} as
\begin{equation}\label{is}
\sum_{k=1}^n\frac{\theta(q^{y_k-1}bz_k)\prod_{j=1}^n\theta
(q^{y_k}z_k/z_j)}
{\theta(q^{N+y_k}bz_k)\prod_{j\neq k}\theta(q^{y_k-y_j}z_k/z_j)}
=\frac{\theta(q^{N+1})\prod_{k=1}^n\theta(q^Nbz_k)}
{\prod_{k=1}^n\theta(q^{N+y_k}bz_k)}.\end{equation}
It follows that
$$R=\sum_{\substack{y_1+\dotsm+y_n=N+1\\y_1,\dots,y_n\geq 0}}
\frac{\De(zq^y)}{\De(z)}
\prod_{k=1}^n\frac{\prod_{j=1}^{n+1}(a_jz_k)_{y_k}}
{(bz_k)_{y_k}\prod_{j=1}^n(qz_k/z_j)_{y_k}}, $$
which completes the proof of the theorem.
\end{proof}

For some
 purposes it is convenient to rewrite  Theorem \ref{taj}  in a way that
hides much of its symmetry but is closer to how the one-variable
Jackson sum is normally written. 
Namely,  we replace $n$ in  Theorem \ref{taj} by $n+1$, eliminate
$y_{n+1}=N-(y_1+\dots+y_n)$ from the summation and write
$z_{n+1}=a^{-1}q^{-N}$. Then the factor
${\De(zq^y)}/{\De(z)}$ is replaced by
$$\frac{\De(zq^y)}{\De(z)}\prod_{k=1}^n\left( q^{N-|y|}
\frac{\theta(z_kq^{y_k}/z_{n+1}q^{N-|y|})}{\theta(z_k/z_{n+1})}\right)
=\frac{\De(zq^y)}{\De(z)}\,q^{n(N-|y|)}\prod_{k=1}^n
\frac{\theta(aq^{y_k+|y|}z_k)}{\theta(aq^Nz_k)}. $$
After  changing  variables $a_j\mapsto b_j$, $b\mapsto aq/c$ 
and using  \eqref{ep} repeatedly, one  obtains the following result.

\begin{corollary}\label{aaj}
If $a^2q^{1+N}=b_1\dotsm b_{n+2}cz_1\dotsm z_n$, the following identity holds:
\begin{multline*}
\sum_{\substack{y_1+\dots+y_n\leq N\\y_1,\dots, y_n\geq 0}}
\Bigg(\frac{\De(zq^y)}{\De(z)}\prod_{k=1}^n\frac{\theta(az_kq^{y_k+|y|})}
{\theta(az_k)}\frac{(q^{-N},c)_{|y|}
\prod_{j=1}^n(az_j)_{|y|}}{\prod_{j=1}^{n+2}(aq/b_j)_{|y|}}\,q^{|y|}\\
\times\prod_{k=1}^n\frac{\prod_{j=1}^{n+2}(b_jz_k)_{y_k}}{(aq^{1+N}z_k,
aqz_k/c)_{y_k}\prod_{j=1}^n(qz_k/z_j)_{y_k}}\Bigg)
=c^N\prod_{k=1}^n\frac{(aqz_k)_N}{(aqz_k/c)_N}\prod_{k=1}^{n+2}
\frac{(aq/cb_k)_N}{(aq/b_k)_N}.
\end{multline*}
\end{corollary}

By a continuation argument, it is possible to convert  Theorem \ref{taj} to 
a sum on a hyper-rectangle.
When $p=0$, this is  \cite[Theorem 6.14]{mr}, given in 
notation more similar to ours as  \cite[Theorem A.5]{mn}. 
The step from the simplex to the hyper-rectangle
 is then based on
polynomial continuation (two polynomials that agree in an infinite
number of points are identical). In the elliptic case, we use
 the quasi-periodicity of the theta function  to pass from an infinite set to
a set having a limit point. (The same argument occurs in \cite{w}.)

 \begin{corollary}\label{cr}
If $a^2q^{1+|m|}=bcde$, the following identity holds:
\begin{equation}\label{ri}\begin{split}&
\sum_{y_1,\dots,y_n=0}^{m_1,\dots,m_n}
\Bigg(\frac{\De(zq^y)}{\De(z)}\prod_{k=1}^n\frac{\theta(az_kq^{y_k+|y|})}
{\theta(az_k)}\frac{(b,c)_{|y|}}{(aq/d,aq/e)_{|y|}}
\prod_{j=1}^n\frac{(az_j)_{|y|}}{(aq^{1+m_j}z_j)_{|y|}}\,q^{|y|}\\
&\qquad
\times\prod_{k=1}^n\frac{(dz_k,ez_k)_{y_k}}{(aqz_k/b,aqz_k/c)_{y_k}}
\prod_{j,k=1}^n\frac{(q^{-m_j}z_k/z_j)_{y_k}}{(qz_k/z_j)_{y_k}}\Bigg)\\
&\quad=\frac{(aq/cd,aq/bd)_{|m|}}{(aq/d,aq/bcd)_{|m|}}
\prod_{k=1}^n\frac{(aqz_k,aqz_k/bc)_{m_k}}{(aqz_k/b,aqz_k/c)_{m_k}}.
\end{split}\end{equation}
\end{corollary}

\begin{proof}
It is straightforward to check that the case $b=q^{-N}$ of Corollary
\ref{cr} is equi\-valent to the case $b_j=q^{-m_j}/z_j$, $1\leq j\leq n$,
of  Corollary \ref{aaj}. Consider the function
$$f(b)=\prod_{k=1}^n(aqz_k/b,aqz_k/c)_{m_k}(L-R),$$
where $L$ and $R$ denote the left- and right-hand sides of \eqref{ri},
respectively,
and where we view $c=a^2q^{1+|m|}/bde$  as depending on $b$ while the
other parameters are fixed. Then $f(b)$ is analytic for $b\neq 0$,
and zero for $b=q^{-N}$, $N\in\mathbb N$.
Since, by \eqref{qp}, we have that in general
$$(\la pb,\la c/p)_k=(c/pb)^k(\la b,\la c)_k, $$
$f$ is quasi-periodic in the sense that
$$f(pb)=(c/pb)^{|m|}f(b). $$
It follows that the points $b=p^kq^{-l}$, $k\in\mathbb Z$, $l\in\mathbb N$
are zeroes of $f$. Assuming that these points are all distinct, that
is, that $p^k\neq q^l$ for $k,\,l\in\mathbb Z$, they have a limit point
in $\mathbb C\setminus \{0\}$ (indeed, in any annulus of the form
$\{z;\,pr\leq |z|\leq r\}$)
and we may conclude that $f\equiv 0$. 
Since $f$ depends analytically on $p$ and $q$, this is true also
if $p^k=q^l$, as long as both sides of \eqref{ri} are well-defined.
\end{proof}

\section{An elliptic $D_n$ Jackson summation}

In this section we give a Jackson-type summation formula
for an elliptic hypergeometric series on the root system $D_n$.
In the case $p=0$, it is due to Schlosser \cite[Theorem 5.17]{s}.
As is discussed below, an essentially equivalent identity was
independently found (still for $p=0$) by Bhatnagar \cite{b2}.

\begin{theorem}\label{tdj} The following identity holds:
\begin{equation}\label{dj}\begin{split}&
\sum_{\substack{y_1+\dots+y_n=N\\y_1,\dots, y_n\geq 0}}
\frac{\De(zq^y)}{\De(z)}\prod_{1\leq j<k\leq n}\frac{1}{(z_jz_k)_{y_j+y_k}}
\prod_{k=1}^n\frac{q^{\binom{y_k}2}z_k^{y_k}\prod_{j=1}^{n-1}
(z_ka_j,z_k/a_j)_{y_k}}{(bz_k,q^{1-N}z_k/b)_{y_k}\prod_{j=1}^n(qz_k/z_j)_{y_k}}
\\
&\quad=(-q^{N-1}b)^N\frac{\prod_{k=1}^{n-1}(ba_k,b/a_k)_N}
{(q)_N\prod_{k=1}^n(bz_k,b/z_k)_N}.
\end{split}\end{equation}
\end{theorem}

\begin{proof}
The proof is by induction on $N$, 
similar to that of Theorem \ref{taj}.

One easily checks that the case $N=1$ of 
\eqref{dj} is equivalent to \eqref{pfd}.
Assume that \eqref{dj} holds for a fixed $N$.  Writing
$R$ for the right-hand side of \eqref{dj} when $N$ is replaced by $N+1$, we  
expand $R$ using first the induction hypothesis  and then 
the case $N=1$ of \eqref{dj} with 
$z_k$ replaced with $z_kq^{y_k}$ and $b$ with $q^Nb$. This gives
\begin{equation*}\begin{split}R&=
(-q^{N}b)^{N+1}\frac{\prod_{k=1}^{n-1}(ba_k,b/a_k)_{N+1}}
{(q)_{N+1}\prod_{k=1}^n(bz_k,b/z_k)_{N+1}}\\
&=-q^{2N}b\,\frac{\prod_{k=1}^{n-1}\theta(q^Nba_k,q^Nb/a_k)}{\theta(q^{N+1})
\prod_{k=1}^n\theta (q^Nbz_k,q^Nb/z_k)}
\,(-q^{N-1}b)^N\frac{\prod_{k=1}^{n-1}(ba_k,b/a_k)_N}
{(q)_N\prod_{k=1}^n(bz_k,b/z_k)_N}\\
&=-q^{2N}b\,\frac{\prod_{k=1}^{n-1}\theta(q^Nba_k,q^Nb/a_k)}{\theta(q^{N+1})
\prod_{k=1}^n\theta (q^Nbz_k,q^Nb/z_k)}\\
&\quad\times\sum_{\substack{y_1+\dots+y_n=N\\y_1,\dots, y_n\geq 0}}
\frac{\De(zq^y)}{\De(z)}\prod_{1\leq j<k\leq n}\frac{1}{(z_jz_k)_{y_j+y_k}}
\prod_{k=1}^n\frac{q^{\binom{y_k}2}z_k^{y_k}\prod_{j=1}^{n-1}
(z_ka_j,z_k/a_j)_{y_k}}{(bz_k,q^{1-N}z_k/b)_{y_k}
\prod_{j=1}^n(qz_k/z_j)_{y_k}}\\
&=q^{N}\frac{\theta (q)}{\theta(q^{N+1})}
\sum_{\substack{y_1+\dots+y_n=N\\y_1,\dots, y_n\geq 0}}\Bigg(
\frac{\De(zq^y)}{\De(z)}\prod_{1\leq j<k\leq n}\frac{1}{(z_jz_k)_{y_j+y_k}}
\\
&\quad\times\prod_{k=1}^n
\frac{\theta(q^{N+y_k}bz_k,q^{N-y_k}b/z_k)\,q^{\binom{y_k}2}z_k^{y_k}
\prod_{j=1}^{n-1}
(z_ka_j,z_k/a_j)_{y_k}}{\theta (q^Nbz_k,q^Nb/z_k)(bz_k,q^{1-N}z_k/b)_{y_k}
\prod_{j=1}^n(qz_k/z_j)_{y_k}}\\
&\quad\times\sum_{\substack{x_1+\dots+x_n=1\\x_1,\dots, x_n\geq 0}}
\bigg\{\frac{\De(zq^{y+x})}{\De(zq^y)}
\prod_{1\leq j<k\leq n}\frac{1}{(z_jz_kq^{y_j+y_k})_{x_j+x_k}}\\
&\quad\times\prod_{k=1}^n\frac{q^{\binom{x_k}2+y_kx_k}z_k^{x_k}
\prod_{j=1}^{n-1}
(z_kq^{y_k}a_j,z_kq^{y_k}/a_j)_{x_k}}
{(q^{N+y_k}bz_k,q^{y_k-N}z_k/b)_{x_k}
\prod_{j=1}^n(q^{1+y_k-y_j}z_k/z_j)_{x_k}}\bigg\}\Bigg).
\end{split}\end{equation*}
Next we replace $y$ by $y-x$ and perform the same  simplifications
as in the proof of Theorem \ref{taj}. We also write
$$\frac{\theta(q^{N-y_k+x_k}b/z_k)}{\theta(q^Nb/z_k)}\frac{1}
{(q^{1-N}z_k/b)_{y_k-x_k}(q^{y_k-x_k-N}z_k/b)_{x_k}}=
\frac{q^{x_k-y_k}}{(q^{-N}z_k/b)_{y_k}}.$$
This gives
\begin{equation*}\begin{split}R&=\frac{1}{\theta(q^{N+1})}
\sum_{\substack{y_1+\dots+y_n=N+1\\y_1,\dots, y_n\geq 0}}\Bigg(
\frac{\De(zq^y)}{\De(z)}\prod_{1\leq j<k\leq n}\frac{1}{(z_jz_k)_{y_j+y_k}}
\\
&\quad\times\prod_{k=1}^n\frac{\theta(q^{N+y_k}bz_k)\,q^{\binom{y_k}2}z_k^{y_k}
\prod_{j=1}^{n-1}
(z_ka_j,z_k/a_j)_{y_k}}{\theta (q^Nbz_k)(bz_k,q^{-N}z_k/b)_{y_k}
\prod_{j=1}^n(qz_k/z_j)_{y_k}}\\
&\quad\times\sum_{\substack{x_1+\dots+x_n=1\\x_1,\dots, x_n\geq 0}}
\prod_{k=1}^n\frac{(q^{y_k-1}bz_k)_{x_k}\prod_{j=1}^n(q^{y_k}z_k/z_j)_{x_k}}
{(q^{N+y_k}bz_k)_{x_k}\prod_{j\neq k}(q^{y_k-y_j}z_k/z_j)_{x_k}}\Bigg).
\end{split}\end{equation*}
Here, the inner sum is identical to the one in \eqref{re} and thus
equals the right-hand side of \eqref{is}. It follows that
$$R=\sum_{\substack{y_1+\dots+y_n=N+1\\y_1,\dots, y_n\geq 0}}
\frac{\De(zq^y)}{\De(z)}\prod_{1\leq j<k\leq n}\frac{1}{(z_jz_k)_{y_j+y_k}}
\prod_{k=1}^n\frac{q^{\binom{y_k}2}z_k^{y_k}\prod_{j=1}^{n-1}
(z_ka_j,z_k/a_j)_{y_k}}{(bz_k,q^{-N}z_k/b)_{y_k}
\prod_{j=1}^n(qz_k/z_j)_{y_k}},$$
which completes the proof.
\end{proof}

Next we rewrite Theorem \ref{tdj} in analogy with
 Corollary \ref{aaj}. This is, up to a  change of variables,
 the form in which the case $p=0$ of 
Theorem \ref{tdj} is  given in \cite{s}.

\begin{corollary}
The following identity holds:
\begin{multline*}
\sum_{\substack{y_1+\dots+y_n\leq N\\y_1,\dots, y_n\geq 0}}
\Bigg(\frac{\De(zq^y)}{\De(z)}\prod_{k=1}^n\frac{\theta(az_kq^{y_k+|y|})}
{\theta(az_k)}\prod_{1\leq j<k\leq n}\frac{1}{(z_jz_k)_{y_j+y_k}}
\prod_{j,k=1}^n\frac{(z_kb_j,z_k/b_j)_{y_k}}{(qz_k/z_j)_{y_k}}
\\
\begin{split}&\quad\times\prod_{j=1}^n\frac{(az_j)_{|y|}(aq/z_j)_{|y|-y_j}}
{(aqb_j,aq/b_j)_{|y|}}\frac{(q^{-N},c,a^2q^{N+1}/c)_{|y|}}
{\prod_{k=1}^n(aqz_k/c,q^{-N}cz_k/a,q^{N+1}az_k)_{y_k}}\,q^{|y|}\Bigg)\\
&=\prod_{k=1}^n\frac{(aqz_k,aq/z_k,aqb_k/c,aq/b_kc)_N}
{(aqz_k/c,aq/z_kc,aqb_k,aq/b_k)_N}.
\end{split}\end{multline*}
\end{corollary}

With exactly the same proof as for Corollary \ref{cr}, 
one may obtain an accompanying
identity for a sum supported on a hyper-rectangle.
When $p=0$, this identity was found independently by Schlosser
\cite[Theorem 5.14]{s} and Bhatnagar \cite[Theorem~7]{b2}.

\begin{corollary}\label{drs}
If $a^2q=bcd$, the following identity holds:
\begin{equation*}\begin{split}&
\sum_{y_1,\dots,y_n=0}^{m_1,\dots,m_n}
\Bigg(\frac{\De(zq^y)}{\De(z)}\prod_{k=1}^n\frac{\theta(az_kq^{y_k+|y|})}
{\theta(az_k)}\prod_{1\leq j<k\leq n}\frac{1}{(z_jz_k)_{y_j+y_k}}
\prod_{j,k=1}^n\frac{(q^{-m_j}z_k/z_j,q^{m_j}z_jz_k)_{y_k}}{(qz_k/z_j)_{y_k}}
\\
&\quad\quad\times\prod_{j=1}^n\frac{(az_j)_{|y|}(aq/z_j)_{|y|-y_j}}
{(aq^{1+m_j}z_j,aq^{1-m_j}/z_j)_{|y|}}\frac{(b,c,d)_{|y|}}
{\prod_{k=1}^n(aqz_k/b,aqz_k/c,aqz_k/d)_{y_k}}\,q^{|y|}\Bigg)\\
&\quad=\prod_{k=1}^n\frac{(aqz_k,bz_k/a,cz_k/a,dz_k/a)_{m_k}}
{(z_k/a,aqz_k/b,aqz_k/c,aqz_k/d)_{m_k}}.
\end{split}\end{equation*}
\end{corollary}

 Corollary \ref{drs} has been written in a form
 similar to how its special case $p=0$ is given in \cite{s}.
In \cite{b2}, it is written down in a different way,
which is obtained by reversing the order of summation. We give this version
of  Corollary \ref{drs} also in the elliptic case, since it is
useful for applications, in particular for obtaining 
the Bailey-type transformations of Section \ref{bs}.

\begin{corollary}\label{drb}
If $a^2q^{1+|m|}=bcde$, the following identity holds:
\begin{equation*}\begin{split}&
\sum_{y_1,\dots,y_n=0}^{m_1,\dots,m_n}
\Bigg(\frac{\De(zq^y)}{\De(z)}\prod_{k=1}^n\frac{\theta(az_kq^{y_k+|y|})}
{\theta(az_k)}\frac{\prod_{1\leq j<k\leq n}(aqz_jz_k/e)_{y_j+y_k}}
{\prod_{j,k=1}^n(aqz_jz_k/e)_{y_k}}
\prod_{j,k=1}^n\frac{(q^{-m_j}z_k/z_j)_{y_k}}{(qz_k/z_j)_{y_k}}
\\
&\quad\quad\times\prod_{j=1}^n\frac{(az_j,e/z_j)_{|y|}}
{(aq^{1+m_j}z_j)_{|y|}(e/z_j)_{|y|-y_j}}\frac
{\prod_{k=1}^n(bz_k,cz_k,dz_k)_{y_k}}{(aq/b,aq/c,aq/d)_{|y|}}\,q^{|y|}\Bigg)\\
&=\frac{\prod_{1\leq j<k\leq n}(aqz_jz_k/e)_{m_j+m_k}}
{\prod_{j,k=1}^n(aqz_jz_k/e)_{m_k}}
\frac{\prod_{k=1}^n(aqz_k,aqz_k/be,aqz_k/ce,aq^{1+|m|-m_k}/bcz_k)_{m_k}}
{(aq/b,aq/c,aq/bce)_{|m|}}.
\end{split}\end{equation*}
\end{corollary}

Note that using \eqref{ip} we may write
$$\frac{\prod_{k=1}^n(aq^{1+|m|-m_k}/bcz_k)_{m_k}}{(aq/bce)_{|m|}}
=
\frac{e^{|m|}q^{\sum_{j<k}m_jm_k}}{z_1^{m_1}\dotsm z_n^{m_n}}
\frac{\prod_{k=1}^n(aqz_k/de)_{m_k}}{(aq/d)_{|m|}},$$
which makes the right-hand side appear more symmetric.

\section{An elliptic  $C_n$ Jackson summation}

Next we give an elliptic Jackson summation on the root system $C_n$.
In the case $p=0$ it was found independently by
Denis and Gustafson \cite[Theorem 4.1]{dg} and by Milne and Lilly 
\cite[Theorem 6.13]{ml}. The general case was stated by
 van Diejen and Spiridonov \cite{ds2}, who showed that it follows
from an elliptic Selberg integral conjectured in \cite{ds}.
They also used modular forms to prove that the Taylor expansion in
$\log q$ of both sides agree up to order $10$
(more precisely, this follows from the non-existence of cusp forms with 
weight less than $12$).

\begin{theorem}\label{cjt}
If $a^2q^{1+|m|}=bcde$, one has the identity
\begin{equation*}\begin{split}&
\sum_{y_1,\dots,y_n=0}^{m_1,\dots,m_n}
\Bigg(\frac{\De(zq^y)}{\De(z)}\prod_{1\leq j\leq k\leq n}
\frac{\theta(az_jz_kq^{y_j+y_k})}{\theta(az_jz_k)}
\prod_{j,k=1}^n\frac{(q^{-m_j}z_k/z_j,az_jz_k)_{y_k}}{(qz_k/z_j,
aq^{1+m_j}z_jz_k)_{y_k}}
\\
&\quad\quad\times\prod_{k=1}^n\frac
{(bz_k,cz_k,dz_k,ez_k)_{y_k}}{(aqz_k/b,aqz_k/c,aqz_k/d,aqz_k/e)_{y_k}}
\,q^{|y|}\Bigg)\\
&=
\frac{\prod_{j,k=1}^n(aqz_jz_k)_{m_k}}
{\prod_{1\leq j<k\leq n}(aqz_jz_k)_{m_j+m_k}}
\frac{(aq/bc,aq/bd,aq/cd)_{|m|}}
{\prod_{k=1}^n(aqz_k/b,aqz_k/c,aqz_k/d,q^{-m_k}e/az_k)_{m_k}}.
\end{split}\end{equation*}
\end{theorem}

The method that we have used above for the root systems $A_n$ and $D_n$
fails for this $C_n$ identity. The reason is that the sum 
 does not involve any
factors of the form $(b)_{|y|}$, so we cannot  convert
it to a sum on a simplex $|y|\leq N$. 
However, the absence of such
factors allows one to employ a different method, namely,
determinant evaluation.

The idea to obtain multivariable hypergeometric summations from
determinant evaluations is due to Gustafson and Krattenthaler \cite{gk}
and was further developped by Schlosser \cite{s2}. Previously, no close
relation has been known between sums coming from determinant evaluations
and the Gustafson--Milne-type  sums that are our main concern here. It 
may therefore seem surprising that we will  obtain
Theorem \ref{cjt} as a special case of a multivariable Jackson sum proved by 
Warnaar \cite{w} using determinant evaluation 
(his ``first'' summation mentioned in the introduction). 
In particular, it follows that the 
 Denis--Gustafson--Milne--Lilly sum is a special case of 
the case $p=0$ of Warnaar's identity, which is due to Schlosser
\cite[Theorem 4.2]{s2}. 

Warnaar's first Jackson summation may be written as
\begin{equation}\label{wj}\begin{split}
&\sum_{y_1,\dots,y_n=0}^N\Bigg(\frac{\Delta(zq^y)}
{\Delta(z)}
\prod_{1\leq j\leq k\leq n} \frac{\theta(az_jz_kq^{y_j+y_k})}
{\theta(az_jz_k)} \\
&\quad\times\prod_{k=1}^n
\frac{(az_k^2,bz_k,cz_k,dz_k,ez_k,q^{-N})_{y_k}}
{(q,aqz_k/b,aqz_k/c,aqz_k/d,aqz_k/e,aq^{N+1}z_k^2)_{y_k}}\,q^{|y|}\Bigg)\\
&=\prod_{1\leq j< k\leq n} \frac{\theta(aq^{N}z_jz_k)}
{\theta(az_jz_k)}
\prod_{k=1}^n\frac{(aqz_k^2,aq^{2-k}/bc,aq^{2-k}/bd,aq^{2-k}/cd)_N}
{(aqz_k/b,aqz_k/c,aqz_k/d,q^{-N}e/az_k)_N},
\end{split}\end{equation}
where $a^2q^{2+N-n}=bcde$. In \cite{r}, we used the case $N=1$ of
this identity to prove Warnaar's conjectured ``second'' 
 elliptic Jackson summation 
(different in nature both from \eqref{wj} and from the results
of the present paper). 
We also remarked that the case $N=1$ of \eqref{wj} is equivalent to
the case $m_j\equiv 1$ of Theorem \ref{cjt}. We will prove 
Theorem \ref{cjt} by combining this observation with the fact 
that  the \emph{general} 
case of Theorem \ref{cjt} is equivalent to
its special case when  $m_j\equiv 1$. This may seem strange but is
actually a common phenomenon: for instance, the summation formulas of
 Corollary \ref{cr} and Corollary \ref{drs}  are also
equivalent to their special case $m_j\equiv 1$.  In a somewhat different
context, this  phenomenon was exploited in 
the proof of  Theorem 3.1 in \cite{r2}.

\begin{proof}[Proof of Theorem \emph{\ref{cjt}}]
We first verify the claim made in \cite{r} that the case $N=1$ 
of \eqref{wj} is equivalent to the case $m_j\equiv 1$ of Theorem \ref{cjt}.

If we let $N=1$ in \eqref{wj} and use \eqref{ip} we obtain
\begin{equation}\label{wj1}\begin{split}
&\sum_{y_1,\dots,y_n=0}^1\Bigg(\frac{\Delta(zq^y)}
{\Delta(z)}
\prod_{1\leq j< k\leq n} \frac{\theta(az_jz_kq^{y_j+y_k})}
{\theta(az_jz_k)} \\
&\quad\times\prod_{k=1}^n\frac{(bz_k,cz_k,dz_k,ez_k)_{y_k}}
{(aqz_k/b,aqz_k/c,aqz_k/d,aqz_k/e)_{y_k}}(-1)^{|y|}\Bigg)\\
&=\frac{\prod_{1\leq j\leq k\leq n}\theta(aqz_jz_k)}
{\prod_{1\leq j< k\leq n}\theta(az_jz_k)}
\prod_{k=1}^n\frac{\theta(aq^{2-k}/bc,aq^{2-k}/bd,aq^{2-k}/cd)}
{\theta(aqz_k/b,aqz_k/c,aqz_k/d,q^{-1}e/az_k)},
\end{split}\end{equation}
where $a^2q^{3-n}=bcde$. 
Next we put  $m_j\equiv 1$ in
 Theorem \ref{cjt} and rewrite the left-hand side
using the two identities
$$\frac{\De(zq^y)}{\De(z)}\prod_{j,k=1}^n\frac{(q^{-1}z_k/z_j)_{y_k}}
{(qz_k/z_j)_{y_k}}=(-1)^{|y|}q^{-|y|}\frac{\De(zq^{-y})}{\De(z)},$$
$$\prod_{1\leq j\leq k\leq n}
\frac{\theta(az_jz_kq^{y_j+y_k})}{\theta(az_jz_k)}
\prod_{j,k=1}^n\frac{(az_jz_k)_{y_k}}{(aq^{2}z_jz_k)_{y_k}}
=\prod_{1\leq j< k\leq n}
\frac{\theta(az_jz_kq^{2-y_j-y_k})}{\theta(az_jz_kq^2)}.
$$
These are easily verified by writing
$$\prod_{j,k=1}^n a_{jk}=\prod_{k=1}^na_{kk}
\prod_{1\leq j< k\leq n}a_{jk}a_{kj} $$
and then considering the four cases $y_j,y_k=0,1$ separately.
Similarly, on the right-hand side we write
$$\frac{\prod_{j,k=1}^n(aqz_jz_k)_{1}}
{\prod_{1\leq j<k\leq n}(aqz_jz_k)_{2}}=
\frac{\prod_{1\leq j\leq k\leq n}\theta(aqz_jz_k)}
{\prod_{1\leq j< k\leq n}\theta(aq^2z_jz_k)}.$$
 This gives
\begin{equation}\label{cj1}\begin{split}
&\sum_{y_1,\dots,y_n=0}^1\Bigg(\frac{\Delta(zq^{-y})}
{\Delta(z)}
\prod_{1\leq j< k\leq n} \frac{\theta(az_jz_kq^{2-y_j-y_k})}
{\theta(az_jz_kq^2)} \\
&\quad\times\prod_{k=1}^n\frac{(bz_k,cz_k,dz_k,ez_k)_{y_k}}
{(aqz_k/b,aqz_k/c,aqz_k/d,aqz_k/e)_{y_k}}(-1)^{|y|}\Bigg)\\
&=\frac{\prod_{1\leq j\leq k\leq n}\theta(aqz_jz_k)}
{\prod_{1\leq j< k\leq n}\theta(aq^2z_jz_k)}
\prod_{k=1}^n
\frac{\theta(aq^k/bc,aq^k/bd,aq^k/cd)}
{\theta(aqz_k/b,aqz_k/c,aqz_k/d,q^{-1}e/az_k)},       
\end{split}\end{equation}
where $a^2q^{1+n}=bcde$.  After the change of variables
 $a\mapsto aq^2$, $q\mapsto q^{-1}$, this
is \eqref{wj1}. (Here we use that, since  $y_k\in\{0,1\}$,
the elliptic Pochhammer symbols  occurring in
\eqref{cj1} do not depend on $q$ for their definition.)

Next we prove that \eqref{cj1} implies the general case 
of Theorem \ref{cjt}. For this we replace $n$ in \eqref{cj1} by
$|m|=m_1+\dots+m_n$ and denote the parameters  $z_j$ in \eqref{cj1} by $w_j$.
We choose these parameters as (compare \cite{r2}, where this
type of choice turned up naturally)
\begin{equation}\label{wz}
(w_1,\dots,w_{|m|})=(z_1,qz_1,\dots,q^{m_1-1}z_1,\dots,
z_n,qz_n,\dots,q^{m_n-1}z_n).\end{equation}
Then the factor $\De(wq^{-y})$ vanishes unless $y$ is of the form
\begin{equation}\label{xy}
y=(\underbrace{\overbrace{1,\dots,1}^{x_1},0,\dots,0}_{m_1},
\dots,\underbrace{\overbrace{1,\dots,1}^{x_n},0,\dots,0}_{m_n}),
\qquad 0\leq x_j\leq m_j. \end{equation}
We claim that if we rewrite \eqref{cj1} using the $x_j$ 
as summation variables, we recover Theorem \ref{cjt}.

The single products in \eqref{cj1} are easily handled using the obvious
identities
$$\prod_{k=1}^{|m|}(bw_k)_{y_k}=\prod_{k=1}^n(bz_k)_{x_k},\qquad
\prod_{k=1}^{|m|}\theta(bw_k)=\prod_{k=1}^n(bz_k)_{m_k}. $$ 
As for the double products, we will prove that, for parameters 
related by \eqref{wz} and \eqref{xy},
\begin{equation}\label{d1}
\prod_{1\leq j< k\leq |m|} \frac{\theta(aw_jw_kq^{2-y_j-y_k})}
{\theta(aw_jw_kq^2)}=\prod_{1\leq j\leq k\leq n}
\frac{\theta(az_jz_kq^{x_j+x_k})}{\theta(az_jz_k)}
\prod_{j,k=1}^n\frac{(az_jz_k)_{x_k}}{(aq^{1+m_j}z_jz_k)_{x_k}},\end{equation}
\begin{equation}\label{d2}\frac{\Delta(wq^{-y})}
{\Delta(w)}=(-1)^{|x|}q^{|x|}\frac{\De(zq^x)}{\De(z)}
\prod_{j,k=1}^n\frac{(q^{-m_j}z_k/z_j)_{x_k}}{(qz_k/z_j)_{x_k}},\end{equation}
\begin{equation}\label{d3}\frac{\prod_{1\leq j\leq k\leq |m|}\theta(aqw_jw_k)}
{\prod_{1\leq j< k\leq |m|}\theta(aq^2w_jw_k)}=
\frac{\prod_{j,k=1}^n(aqz_jz_k)_{m_k}}
{\prod_{1\leq j<k\leq n}(aqz_jz_k)_{m_j+m_k}}.\end{equation}
Assuming that these identities have been proved it is easily checked that
our claim is correct, so that  Theorem \ref{cjt} follows.

The left-hand sides of \eqref{d1} and \eqref{d2} are both of the  form 
$$\prod_{1\leq j<k\leq |m|}\frac{f(w_jq^{-y_j},w_kq^{-y_k})}{f(w_j,w_k)}.$$
We decompose these products as $ABC$, where
$$A=\prod_{j,k=1}^n\prod_{u=1}^{x_k}\prod_{t=x_j+1}^{m_j}
\frac{f(z_jq^{t-1},z_kq^{u-2})}{f(z_jq^{t-1},z_kq^{u-1})},$$
$$B=\prod_{k=1}^n\prod_{1\leq t<u\leq x_k}
\frac{f(z_kq^{t-2},z_kq^{u-2})}{f(z_kq^{t-1},z_kq^{u-1})}, $$
$$C=\prod_{1\leq j<k\leq n}\prod_{u=1}^{x_k}\prod_{t=1}^{x_j}
\frac{f(z_jq^{t-2},z_kq^{u-2})}{f(z_jq^{t-1},z_kq^{u-1})}. $$
Here we use the symmetry and antisymmetry, respectively, of $f$ 
to collect the factors with $y_j\neq y_k$ as in $A$.

In the case of \eqref{d1}, we have
\begin{equation}\label{ap1}\begin{split}
A&=\prod_{j,k=1}^n\prod_{u=1}^{x_k}\prod_{t=x_j+1}^{m_j}
\frac{\theta(az_jz_kq^{t+u-1})}{\theta(az_jz_kq^{t+u})}
=\prod_{j,k=1}^n\prod_{u=1}^{x_k}
\frac{\theta(az_jz_kq^{x_j+u})}{\theta(az_jz_kq^{m_j+u})}\\
&=\prod_{j,k=1}^n
\frac{(az_jz_kq^{1+x_j})_{x_k}}{(az_jz_kq^{1+m_j})_{x_k}}.
\end{split}\end{equation}
The factor $B$ may be computed using the elementary identity
$$\prod_{1\leq j<k\leq n}\frac{f(j+k-2)}{f(j+k)}=
\prod_{k=1}^n\frac{f(k)}{f(n+k-1)},$$
 which gives
\begin{equation}\label{ap2}B=\prod_{k=1}^n\prod_{1\leq t<u\leq x_k}
\frac{\theta(az_k^2q^{t+u-2})}{\theta(az_k^2q^{t+u})}
=\prod_{k=1}^n\frac{(az_k^2q)_{x_k}}{(az_k^2q^{x_k})_{x_k}}.\end{equation}
Finally, $C$ may be computed similarly as $A$:
\begin{equation}\label{ap3}\begin{split}
C&=\prod_{1\leq j<k\leq n}\prod_{u=1}^{x_k}\prod_{t=1}^{x_j}
\frac{\theta(az_jz_kq^{t+u-2})}{\theta(az_jz_kq^{t+u})}
=\prod_{1\leq j<k\leq n}\prod_{u=1}^{x_k}
\frac{\theta(az_jz_kq^{u-1},az_jz_kq^{u})}{\theta(az_jz_kq^{x_j+u-1},
az_jz_kq^{x_j+u})}\\
&=\prod_{1\leq j<k\leq n}\frac{(az_jz_k,az_jz_kq)_{x_k}}
{(az_jz_kq^{x_j},az_jz_kq^{1+x_j})_{x_k}}.\end{split}\end{equation}
Combining \eqref{ap1}, \eqref{ap2} and \eqref{ap3} gives
\begin{equation*}\begin{split}
\prod_{1\leq j< k\leq |m|} \frac{\theta(aw_jw_kq^{2-y_j-y_k})}
{\theta(aw_jw_kq^2)}&=\prod_{j,k=1}^n
\frac{(az_jz_kq^{1+x_j})_{x_k}}{(az_jz_kq^{1+m_j})_{x_k}}
\prod_{1\leq j\leq k\leq n}\frac{(az_jz_kq)_{x_k}}{(az_jz_kq^{x_j})_{x_k}}\\
&\quad\times\prod_{1\leq j<k\leq n}\frac{(az_jz_k)_{x_k}}
{(az_jz_kq^{1+x_j})_{x_k}}.\end{split}\end{equation*}
To complete the proof of \eqref{d1}, we must show that
\begin{equation*}\begin{split}&\prod_{j,k=1}^n
\frac{(az_jz_kq^{1+x_j})_{x_k}}{(az_jz_k)_{x_k}}
\prod_{1\leq j\leq k\leq n}\frac{(az_jz_kq)_{x_k}}{(az_jz_kq^{x_j})_{x_k}}
\prod_{1\leq j<k\leq n}\frac{(az_jz_k)_{x_k}}
{(az_jz_kq^{1+x_j})_{x_k}}\\
&\qquad=\prod_{1\leq j\leq k\leq n}
\frac{\theta(az_jz_kq^{x_j+x_k})}{\theta(az_jz_k)},\end{split}\end{equation*}
which is easily verified after writing
$$\prod_{j,k=1}^n
\frac{(az_jz_kq^{1+x_j})_{x_k}}{(az_jz_k)_{x_k}}
=\prod_{1\leq j<k\leq n}\frac{(az_jz_kq^{1+x_j})_{x_k}}{(az_jz_k)_{x_k}}
\prod_{1\leq j\leq k\leq n}\frac{(az_jz_kq^{1+x_k})_{x_j}}{(az_jz_k)_{x_j}}.$$

Equation \eqref{d2} may be proved in a  similar way, and we do
not give the details.  The case  $p=0$, which is essentially the same,  
was done as part of the proof of Theorem 3.1 in \cite{r2}. 

To prove \eqref{d3}, we write the left-hand side as
$$\frac{\prod_{1\leq j\leq k\leq |m|}\theta(aqw_jw_k)}
{\prod_{1\leq j<k\leq |m|}\theta(aq^2w_jw_k)}
=\prod_{k=1}^n\frac{\prod_{1\leq t\leq u\leq m_k}\theta(aq^{t+u-1}z_k^2)}
{\prod_{1\leq t< u\leq m_k}\theta(aq^{t+u}z_k^2)}
\prod_{1\leq j<k\leq n}\prod_{u=1}^{m_k}\prod_{t=1}^{m_j}
\frac{\theta(aq^{t+u-1}z_jz_k)}{\theta(aq^{t+u}z_jz_k)}.$$
Here, the first part of the product may be computed using the identity
$$\frac{\prod_{1\leq j\leq k\leq n}f(j+k-1)}
{\prod_{1\leq j< k\leq n}f(j+k)}=\prod_{k=1}^nf(k), $$
giving
\begin{equation}\label{a}
\prod_{k=1}^n\frac{\prod_{1\leq t\leq u\leq m_k}\theta(aq^{t+u-1}z_k^2)}
{\prod_{1\leq t< u\leq m_k}\theta(aq^{t+u}z_k^2)}=\prod_{k=1}^n(aqz_k^2)_{m_k}.
\end{equation}
Finally, the second part is given by
\begin{multline}\label{y}
\prod_{1\leq j<k\leq n}\prod_{u=1}^{m_k}\prod_{t=1}^{m_j}
\frac{\theta(aq^{t+u-1}z_jz_k)}{\theta(aq^{t+u}z_jz_k)}
=\prod_{1\leq j<k\leq n}\prod_{u=1}^{m_k}
\frac{\theta(aq^{u}z_jz_k)}{\theta(aq^{u+m_j}z_jz_k)}\\
=\prod_{1\leq j<k\leq n}\frac{(aqz_jz_k)_{m_k}(aqz_jz_k)_{m_j}}
{(aqz_jz_k)_{m_j+m_k}}=\frac{1}{\prod_{k=1}^n(aqz_k^2)_{m_k}}
\frac{\prod_{j,k=1}^n(aqz_jz_k)_{m_k}}
{\prod_{1\leq j<k\leq n}(aqz_jz_k)_{m_j+m_k}}.\end{multline}
Combining \eqref{a} and \eqref{y} we obtain \eqref{d3}. This completes
the proof of Theorem~\ref{cjt}.
\end{proof}

\section{Elliptic Bailey transformations}
\label{bs}

As applications of our multivariable  elliptic Jackson summations, 
we may obtain a host of multivariable Bailey transformations
for elliptic hypergeometric series on root systems. It should be observed that
the proofs  are almost identical to those for the
case $p=0$ given in \cite{bs,mn}, 
since they mainly involve  manipulation of Pochhammer
symbols using the elementary identities \eqref{ap}, \eqref{ep}, \eqref{ip},
which do not depend on $p$. The only difference is that one must 
replace  polynomial continuation arguments by  appealing to
the quasi-periodicity, as in the proof Corollary \ref{cr}. 
There is therefore no need to give detailed
 proofs of these new elliptic Bailey transformations. Anyway, we sketch
the proof of the first one to indicate to the reader what is involved;
for the remaining ones we provide enough details so that the sufficiently
interested reader should have no trouble checking the computations.

We begin with a Bailey transformation for elliptic hypergeometric series on 
$A_n$. When $p=0$, it is due to Denis and Gustafson \cite[Theorem 3.1]{dg},
who derived it by residue calculus from a multivariable integral
transformation. It was rediscovered by 
Milne and Newcomb \cite[Theorem 3.1]{mn}, who used the method  sketched
below.  A third proof was given
in \cite{r2}. (In \cite{mn}, it is erroneously claimed that the
transformation in \cite{mn} has one more  free parameter than the one in 
\cite{dg}, whereas in fact the two results are  equivalent. 
The reason for this mistake seems to be that in \cite[Theorem 3.1]{mn} one 
may multiply all the parameters $x_i$ by a constant
without changing the result. This causes the authors of \cite{mn} to
overestimate the number of free parameters in their identity.)

\begin{corollary}\label{bt}
Assuming $a^3q^{2+|m|}=bcdefg$ and writing $\la=a^2q/bce$, the following 
identity holds:
\begin{equation*}\begin{split}&\sum_{y_1,\dots,y_n=0}^{m_1,\dots,m_n}
\Bigg(\frac{\De(zq^y)}{\De(z)}\prod_{k=1}^n\frac{\theta(az_kq^{y_k+|y|})}
{\theta(az_k)}\frac{(b,c,d)_{|y|}}{(aq/e,aq/f,aq/g)_{|y|}}
\prod_{j=1}^n\frac{(az_j)_{|y|}}{(aq^{1+m_j}z_j)_{|y|}}\,q^{|y|}\\
&\quad\times\prod_{k=1}^n\frac{(ez_k,fz_k,gz_k)_{y_k}}
{(aqz_k/b,aqz_k/c,aqz_k/d)_{y_k}}
\prod_{j,k=1}^n\frac{(q^{-m_j}z_k/z_j)_{y_k}}{(qz_k/z_j)_{y_k}}\Bigg)\\
&=\left(\frac a \la\right)^{|m|}
\frac{(\la q/f,\la q/g)_{|m|}}{(aq/f,aq/g)_{|m|}}
\prod_{k=1}^n\frac{(aqz_k,\la qz_k/d)_{m_k}}{(\la qz_k,aqz_k/d)_{m_k}}\\
&\quad\times\sum_{y_1,\dots,y_n=0}^{m_1,\dots,m_n}
\Bigg(\frac{\De(zq^y)}{\De(z)}\prod_{k=1}^n\frac{\theta(\la z_k q^{y_k+|y|})}
{\theta(\la z_k)}
\frac{(\la b/a,\la c/a,d)_{|y|}}
{(aq/e,\la q/f,\la q/g)_{|y|}}
\prod_{j=1}^n\frac{(\la z_j)_{|y|}}{(\la q^{1+m_j}z_j)_{|y|}}
\,q^{|y|}\\
&\quad\times
\prod_{k=1}^n\frac{(\la ez_k/a,fz_k,gz_k)_{y_k}}
{(aqz_k/b,aqz_k/c,\la qz_k/d)_{y_k}}
\prod_{j,k=1}^n\frac{(q^{-m_j}z_k/z_j)_{y_k}}{(qz_k/z_j)_{y_k}}\Bigg).
\end{split}\end{equation*}
\end{corollary}

\begin{proof}
Let $W$ denote the sum on the left-hand side. 
We want to expand each term by identifying the factor
$$\frac{(b,c)_{|y|}}{(aq/e)_{|y|}}\prod_{k=1}^n\frac{(ez_k)_{y_k}}
{(aqz_k/b,aqz_k/c)_{y_k}} $$
with a part of the right-hand side of \eqref{ri}. 
We must then replace the parameters $(a,b,c,d,e,m_k)$ in 
Corollary \ref{cr}
with  $(\la, \la b/a, \la c/a, \la e/a, aq^{|y|}, y_k)$. 
This gives
\begin{equation*}\begin{split}W&=\sum_{y_1,\dots,y_n=0}^{m_1,\dots,m_n}
\Bigg(\frac{\De(zq^y)}{\De(z)}\prod_{k=1}^n\frac{\theta(az_kq^{y_k+|y|})}
{\theta(az_k)}\frac{(d,a/\la)_{|y|}}{(aq/f,aq/g)_{|y|}}
\prod_{j=1}^n\frac{(az_j)_{|y|}}{(aq^{1+m_j}z_j)_{|y|}}\,q^{|y|}\\
&\quad\times\prod_{k=1}^n\frac{(fz_k,gz_k)_{y_k}}
{(aqz_k/d,\la qz_k)_{y_k}}
\prod_{j,k=1}^n\frac{(q^{-m_j}z_k/z_j)_{y_k}}{(qz_k/z_j)_{y_k}}\\
&\quad\times
\sum_{x_1,\dots,x_n=0}^{y_1,\dots,y_n}
\bigg\{\frac{\De(zq^x)}{\De(z)}\prod_{k=1}^n\frac{\theta(\la z_kq^{x_k+|x|})}
{\theta(\la z_k)}\frac{(\la b/a,\la c/a)_{|x|}}{(aq/e,q^{1-|y|}\la/a)_{|x|}}
\prod_{j=1}^n\frac{(\la z_j)_{|x|}}{(\la q^{1+y_j}z_j)_{|x|}}\,q^{|x|}\\
&\quad\times\prod_{k=1}^n\frac{(\la e z_k/a,aq^{|y|}z_k)_{x_k}}
{(aqz_k/b,aqz_k/c)_{x_k}}
\prod_{j,k=1}^n\frac{(q^{-y_j}z_k/z_j)_{x_k}}{(qz_k/z_j)_{x_k}}\bigg\}
\Bigg).
\end{split}\end{equation*}
Replacing $y$ by $y+x$ and changing the order of summation gives,
after some elementary manipulations and an application of \eqref{di},
\begin{multline*}
\begin{split}W&=\sum_{x_1,\dots,x_n=0}^{m_1,\dots,m_n}
\Bigg(\frac{\De(zq^x)}{\De(z)}\frac{(\la b/a,\la c/a,d)_{|x|}}
{(aq/e,aq/f,aq/g)_{|x|}}
\prod_{j=1}^n\frac{(aqz_j)_{|x|+x_j}}{(aq^{1+m_j}z_j)_{|x|}
(\la q^{|x|}z_j)_{x_j}}\,q^{|x|}\\
&\quad\times
\left(\frac a\la\right)^{|x|}
\prod_{k=1}^n\frac{(\la ez_k/a,fz_k,gz_k)_{x_k}}
{(aqz_k/b,aqz_k/c,aqz_k/d)_{x_k}}
\prod_{j,k=1}^n\frac{(q^{-m_j}z_k/z_j)_{x_k}}
{(qz_k/z_j)_{x_k}}\\
&\quad\times
\sum_{y_1,\dots,y_n=0}^{m_1-x_1,\dots,m_n-x_n}
\bigg\{\frac{\De(zq^{x+y})}{\De(zq^x)}
\prod_{k=1}^n\frac{\theta(a z_kq^{x_k+|x|+y_k+|y|})}
{\theta(az_kq^{x_k+|x|})}\frac{(dq^{|x|},a/\la)_{|y|}}{(aq^{1+|x|}/f,
aq^{1+|x|}/g)_{|y|}}\,q^{|y|}\end{split}\\
\times\prod_{j=1}^n\frac{(aq^{|x|+x_j} z_j)_{|y|}}
{(a q^{1+|x|+m_j}z_j)_{|y|}}
\prod_{k=1}^n\frac{(fq^{x_k}z_k,gq^{x_k}z_k)_{y_k}}
{(aq^{1+x_k}z_k/d,\la q^{1+|x|+x_k}z_k)_{y_k}}
\prod_{j,k=1}^n\frac{(q^{x_k-m_j}z_k/z_j)_{y_k}}
{(q^{1+x_k-x_j}z_k/z_j)_{y_k}}\bigg\}\Bigg).
\end{multline*}
The condition   $bcdefg=a^3q^{2+|m|}$ means that the inner sum 
 is as in Corollary \ref{cr}, with $(a,b,c,d,e,m_k,z_k)$ replaced by
$(aq^{|x|},dq^{|x|},a/\la,f,g,m_k-x_x,q^{x_k}z_k)$. Using
 \eqref{ep} and the fact that $dfg=a\la q^{1+|m|}$,
the corresponding right-hand side of \eqref{ri} may be written as
$$\left(\frac a \la\right)^{|m|-|x|}\frac{(\la q^{1+|x|}/f,
\la q^{1+|x|}/g)_{|m|-|x|}}{(a q^{1+|x|}/f,a q^{1+|x|}/g)_{|m|-|x|}}
\prod_{k=1}^n\frac{(aq^{1+|x|+x_k}z_k,\la q^{1+x_k}z_k/d)_{m_k-x_k}}
{(\la q^{1+|x|+x_k}z_k,a q^{1+x_k}z_k/d)_{m_k-x_k}}.$$
Some further elementary manipulations  completes the proof.
\end{proof}

There is also a  version of Corollary \ref{bt} for series supported
on a simplex. It can be obtained either
by a slight modification of the proof of  Corollary \ref{bt}, using first
Corollary \ref{cr} and then Corollary \ref{aaj}, instead of using Corollary
\ref{cr} twice, or as a consequence of Corollary \ref{bt}, using a
continuation argument similar to the proof of Corollary \ref{cr}.
When $p=0$, the resulting identity is   \cite[Theorem 3.3]{mn}.

\begin{corollary}
Assuming $a^3q^{2+N}=b_1\dotsm b_{n+2}cdez_1\dotsm z_n$
and writing $\la=a^2q/cde$, the following identity holds:
\begin{equation*}\begin{split}&
\sum_{\substack{y_1+\dots+y_n\leq N\\y_1,\dots, y_n\geq 0}}
\Bigg(\frac{\De(zq^y)}{\De(z)}\prod_{k=1}^n\frac{\theta(az_kq^{y_k+|y|})}
{\theta(az_k)}\frac{(q^{-N},c,d)_{|y|}
\prod_{j=1}^n(az_j)_{|y|}}{(aq/e)_{|y|}
\prod_{j=1}^{n+2}(aq/b_j)_{|y|}}\,q^{|y|}\\
&\quad
\times\prod_{k=1}^n\frac{(ez_k)_{y_k}\prod_{j=1}^{n+2}(b_jz_k)_{y_k}}
{(aq^{1+N}z_k,aqz_k/c,aqz_k/d)_{y_k}\prod_{j=1}^n(qz_k/z_j)_{y_k}}\Bigg)\\
&=\left(\frac a\la\right)^N
\prod_{k=1}^n\frac{(aqz_k)_N\prod_{j=1}^{n+2}(\la q/b_j)_N}
{(\la qz_k)_N\prod_{j=1}^{n+2}(aq/b_j)_N}\\
&\quad\times\sum_{\substack{y_1+\dots+y_n\leq N\\y_1,\dots, y_n\geq 0}}
\Bigg(\frac{\De(zq^y)}{\De(z)}\prod_{k=1}^n\frac{\theta(\la z_kq^{y_k+|y|})}
{\theta(\la z_k)}\frac{(q^{-N},\la c/a,\la d/a)_{|y|}
\prod_{j=1}^n(\la z_j)_{|y|}}{(aq/e)_{|y|}
\prod_{j=1}^{n+2}(\la q/b_j)_{|y|}}\,q^{|y|}\\
&\quad
\times\prod_{k=1}^n\frac{(\la ez_k/a)_{y_k}\prod_{j=1}^{n+2}(b_jz_k)_{y_k}}
{(\la q^{1+N}z_k,aqz_k/c,aqz_k/d)_{y_k}\prod_{j=1}^n(qz_k/z_j)_{y_k}}\Bigg).
\end{split}\end{equation*}
\end{corollary}

In \cite{bs} a  number of $C_n$ and $D_n$ Bailey transformations were 
obtained by judiciously combining $A_n$, $C_n$ and $D_n$ Jackson summations,
similarly as in the proof of Corollary \ref{bt}. 
Starting from  the elliptic
Jackson summations obtained in the present paper, the same  
method yields elliptic $C_n$ and $D_n$ Bailey transformations.
The seven transformations given in \cite{bs} fall into three groups.
We are content with writing down one representative from each group
explicitly. 

Combining a $C_n$ and a $D_n$ Jackson summation, one may prove the
following identity relating 
 an elliptic ${}_{10}W_9$ series on $C_n$ with a 
similar series on $A_n$. For $p=0$ it is Theorem 2.1 of \cite{bs}.

\begin{corollary}\label{cbt}
Assuming that $a^3q^{2+|m|}=bcdefg$ and writing $\la=a^2q/bcd$, the following 
identity holds:
\begin{equation*}\begin{split}&
\sum_{y_1,\dots,y_n=0}^{m_1,\dots,m_n}
\Bigg(\frac{\De(zq^y)}{\De(z)}\prod_{1\leq j\leq k\leq n}
\frac{\theta(az_jz_kq^{y_j+y_k})}{\theta(az_jz_k)}
\prod_{j,k=1}^n\frac{(q^{-m_j}z_k/z_j,az_jz_k)_{y_k}}{(qz_k/z_j,
aq^{1+m_j}z_jz_k)_{y_k}}
\\
&\quad\times\prod_{k=1}^n\frac
{(bz_k,cz_k,dz_k,ez_k,fz_k,gz_k)_{y_k}}
{(aqz_k/b,aqz_k/c,aqz_k/d,aqz_k/e,aqz_k/f,aqz_k/g)_{y_k}}
\,q^{|y|}\Bigg)
\\
&=
\frac{\prod_{j,k=1}^n(aqz_jz_k)_{m_k}}
{\prod_{1\leq j<k\leq n}(aqz_jz_k)_{m_j+m_k}}
\frac{(\la q/e,\la q/f,aq/ef)_{|m|}}
{\prod_{k=1}^n(\la qz_k,aqz_k/e,aqz_k/f,q^{-m_k}g/az_k)_{m_k}}\\
&\quad\times\sum_{y_1,\dots,y_n=0}^{m_1,\dots,m_n}
\Bigg(\frac{\De(zq^y)}{\De(z)}\prod_{k=1}^n\frac{\theta(\la z_k q^{y_k+|y|})}
{\theta(\la z_k)}
\frac{(\la b/a,\la c/a,\la d/a)_{|y|}}
{(\la q/e,\la q/f,\la q/g)_{|y|}}
\prod_{j=1}^n\frac{(\la z_j)_{|y|}}{(\la q^{1+m_j}z_j)_{|y|}}
\,q^{|y|}\\
&\quad\times
\prod_{k=1}^n\frac{(ez_k,fz_k,gz_k)_{y_k}}
{(aqz_k/b,aqz_k/c,aqz_k/d)_{y_k}}
\prod_{j,k=1}^n\frac{(q^{-m_j}z_k/z_j)_{y_k}}{(qz_k/z_j)_{y_k}}\Bigg).
\end{split}\end{equation*}
\end{corollary}

Using  \eqref{ip} one may write
$$\frac{(aq/ef)_{|m|}}{\prod_{k=1}^n(q^{-m_k}g/az_k)_{m_k}}
=\left(\frac a\la\right)^{|m|}q^{-\sum_{j<k}m_jm_k}z_1^{m_1}\dotsm z_n^{m_n}
\frac{(\la q/g)_{|m|}}{\prod_{k=1}^n(aqz_k/g)_{m_k}}, $$
which makes the right-hand side appear more symmetric.
Note also  that, since the left-hand side 
of Corollary \ref{cbt} is 
invariant under interchanging $d$ and $g$, this must be
true also for the right-hand side. This fact is equivalent to
 Corollary \ref{bt}, which may thus alternatively be derived as
a consequence of  Corollary \ref{cbt}.

To prove Corollary  \ref{cbt}, one may use the case
$$(a,b,c,d,z_k,m_k)\mapsto(\la/\sqrt{a},\la b/a,\la c/a,\la d/a,
\sqrt a z_k,y_k)$$
of Corollary \ref{drs} to expand the factor
$$\prod_{k=1}^n\frac
{(bz_k,cz_k,dz_k)_{y_k}}{(aqz_k/b,aqz_k/c,aqz_k/d)_{y_k}} $$
on the left-hand side. After changing the order of summation as
in the proof of Corollary \ref{bt}, the inner sum may be computed 
using Theorem \ref{cjt}. 

Next we give a  $D_n$ Bailey transformation which may be obtained
by using a $D_n$ Jackson summation twice.
 For $p=0$ it is equivalent to  \cite[Theorem 3.1]{bs}.

\begin{corollary}\label{yc}
Assuming $a^3q^2=bcdef$ and writing $\la=a^2q/bcd$, the following identity
holds:
\begin{equation*}\begin{split}&
\sum_{y_1,\dots,y_n=0}^{m_1,\dots,m_n}
\Bigg(\frac{\De(zq^y)}{\De(z)}\prod_{k=1}^n\frac{\theta(az_kq^{y_k+|y|})}
{\theta(az_k)}\prod_{j=1}^n\frac{(az_j,d/z_j)_{|y|}(aq/z_j)_{|y|-y_j}}
{(aq^{1+m_j}z_j,aq^{1-m_j}/z_j)_{|y|}(d/z_j)_{|y|-y_j}}\\
&\quad\quad\times
\prod_{j,k=1}^n\frac{(q^{-m_j}z_k/z_j,q^{m_j}z_jz_k)_{y_k}}{(qz_k/z_j,
aqz_jz_k/d)_{y_k}}\prod_{1\leq j<k\leq n}\frac{(aqz_jz_k/d)_{y_j+y_k}}
{(z_jz_k)_{y_j+y_k}}\\
&\quad\quad\times\prod_{k=1}^n
\frac{(bz_k,cz_k)_{y_k}}
{(aqz_k/e,aqz_k/f)_{y_k}}\frac{(e,f)_{|y|}}{(aq/b,aq/c)_{|y|}}\,q^{|y|}\Bigg)\\
&\quad=\prod_{k=1}^n\frac{(aqz_k,z_k/\la,\la qz_k/e,\la qz_k/f)_{m_k}}
{(\la qz_k,z_k/a,aqz_k/e,aqz_k/f)_{m_k}}\sum_{y_1,\dots,y_n=0}^{m_1,\dots,m_n}
\Bigg(\frac{\De(zq^y)}{\De(z)}\prod_{k=1}^n\frac{\theta(\la z_kq^{y_k+|y|})}
{\theta(\la z_k)}\\
&\quad\quad\times\prod_{j=1}^n\frac{(\la z_j,\la d/az_j)_{|y|}
(\la q/z_j)_{|y|-y_j}}
{(\la q^{1+m_j}z_j,\la q^{1-m_j}/z_j)_{|y|}(\la d/az_j)_{|y|-y_j}}
\prod_{j,k=1}^n\frac{(q^{-m_j}z_k/z_j,q^{m_j}z_jz_k)_{y_k}}{(qz_k/z_j,
aqz_jz_k/d)_{y_k}}\\
&\quad\quad\times\prod_{1\leq j<k\leq n}\frac{(aqz_jz_k/d)_{y_j+y_k}}
{(z_jz_k)_{y_j+y_k}}
\prod_{k=1}^n
\frac{(\la bz_k/a,\la cz_k/a)_{y_k}}
{(\la qz_k/e,\la qz_k/f)_{y_k}}\frac{(e,f)_{|y|}}{(aq/b,aq/c)_{|y|}}
\,q^{|y|}\Bigg).
\end{split}\end{equation*}
\end{corollary}

At a first glance, \cite[Theorem 3.1]{bs} appears to have one more free 
parameter than Corollary \ref{yc}, but, as noted in \cite[Remark 3.3]{bs}, 
one may specialize one parameter in that result without loss of generality.

To prove Corollary \ref{yc} one may use the case
$$(a,b,c,d,e,m_k)\mapsto(\la,\la b/a,\la c/a,aq^{|y|},\la d/a,
 y_k) $$
of  Corollary \ref{drb} to expand the factor 
$$\frac{\prod_{1\leq j<k\leq n}(aqz_jz_k/d)_{y_j+y_k}}{\prod_{j,k=1}^n
(aqz_jz_k/d)_{y_k}}\frac{\prod_{k=1}^n(bz_k,cz_k,dq^{|y|-y_k}/z_k)_{y_k}}
{(aq/b,aq/c)_{|y|}}. $$
Proceeding as in the proof of Corollary \ref{bt} calls for an application 
 of Corollary \ref{drs} in the last step.

By a continuation argument, one may obtain a  companion identity
to Corollary \ref{yc} where the sum is over a simplex $|y|\leq N$.
We do not write it out explicitly; the case $p=0$ is 
Theorem 3.7 of \cite{bs}. 

Another class of $D_n$ Bailey transformations may be obtained by combining
an $A_n$ and a $D_n$ Jackson summation.
An example is the following identity, which is equivalent to Theorem 
3.13 of \cite{bs} when $p=0$.

\begin{corollary}\label{dt1}
Assuming that $a^3q^2=bcdef$ and writing
$\la=a^2q/bcf$, the following identity holds:
\begin{multline*}
\sum_{y_1,\dots,y_n=0}^{m_1,\dots,m_n}
\Bigg(\frac{\De(zq^y)}{\De(z)}\prod_{k=1}^n\frac{\theta(az_kq^{y_k+|y|})}
{\theta(az_k)}\prod_{1\leq j<k\leq n}\frac{1}{(z_jz_k)_{y_j+y_k}}
\prod_{j,k=1}^n\frac{(q^{-m_j}z_k/z_j,q^{m_j}z_jz_k)_{y_k}}{(qz_k/z_j)_{y_k}}
\\
\times\prod_{j=1}^n\frac{(az_j)_{|y|}(aq/z_j)_{|y|-y_j}}
{(aq^{1+m_j}z_j,aq^{1-m_j}/z_j)_{|y|}}\frac{(b,c,d,e)_{|y|}}{(aq/f)_{|y|}}
\prod_{k=1}^n\frac{(fz_k)_{y_k}}{(aqz_k/b,aqz_k/c,aqz_k/d,aqz_k/e)_{y_k}}
\,q^{|y|}\Bigg)\\
\begin{split}&=\prod_{k=1}^n\frac{(aqz_k,z_k/\la,\la qz_k/d,\la qz_k/e)_{m_k}}
{(\la qz_k,z_k/a,aqz_k/d,aqz_k/e)_{m_k}}\sum_{y_1,\dots,y_n=0}^{m_1,\dots,m_n}
\Bigg(\frac{\De(zq^y)}{\De(z)}\prod_{k=1}^n\frac{\theta(\la z_kq^{y_k+|y|})}
{\theta(\la z_k)}\\
&\quad\times\prod_{1\leq j<k\leq n}\frac{1}{(z_jz_k)_{y_j+y_k}}
\prod_{j,k=1}^n\frac{(q^{-m_j}z_k/z_j,q^{m_j}z_jz_k)_{y_k}}{(qz_k/z_j)_{y_k}}
\prod_{j=1}^n\frac{(\la z_j)_{|y|}
(\la q/z_j)_{|y|-y_j}}
{(\la q^{1+m_j}z_j,\la q^{1-m_j}/z_j)_{|y|}}
\\
&\quad\times\frac{(\la b/a,\la c/a,d,e)_{|y|}}{(aq/f)_{|y|}}
\prod_{k=1}^n\frac{(\la fz_k/a)_{y_k}}{(aqz_k/b,aqz_k/c,\la qz_k/d,
\la qz_k/e)_{y_k}}\,q^{|y|}\Bigg).
\end{split}\end{multline*}
\end{corollary}

To prove this one may  use the case $(a,b,c,d,e,m_k)\mapsto
(\la,\la b/a,\la c/a,\la f/a,aq^{|y|},y_k)$
 of Corollary \ref{cr} to expand the factor
$$ \frac{(b,c)_{|y|}}{(aq/f)_{|y|}}
\prod_{k=1}^n\frac{(fz_k)_{y_k}}{(aqz_k/b,aqz_k/c)_{y_k}}$$
on the left-hand side. The same method as before leads to a sum
that is computed by Corollary \ref{drs}.

In \cite{bs}, three companion identities to the case $p=0$ of
Corollary \ref{dt1} are given. One of these \cite[Theorem 3.9]{bs}
is the equivalent identity obtained by reversing the order of summation, 
the other two \cite[Theorem 3.11 and Theorem 3.16]{bs} 
(which are not equivalent) 
are continuations of the first two to the simplex. Again, these
three identities are easily extended  to the elliptic case, but
we do not write them out explicitly.

\end{document}